\PassOptionsToPackage{dvipsnames}{xcolor} 
\documentclass[11pt]{article}
\usepackage{amssymb,amsmath,amsthm,multirow,color,xcolor,cancel,bbm,tikz,mathrsfs,hyperref,authblk,stmaryrd}
\usepackage[textheight=24cm,textwidth=16cm]{geometry}
\tikzstyle{main node}=[draw,circle,inner sep=1,outer sep=2,thick,minimum size=12pt]
\usetikzlibrary{decorations.pathreplacing,calligraphy}

\allowdisplaybreaks 

\hypersetup{
    colorlinks=true,
    linkcolor=blue,
    citecolor=cyan,
    filecolor=magenta,      
    urlcolor=green,
    pdftitle={Interaction graphs of isomorphic automata networks I},
    }

\newtheorem{definition}{Definition}

\newtheorem{lemma}{Lemma}
\newtheorem{theorem}{Theorem}

\newcommand{\Q}[1]{\llbracket #1\rrbracket}
\newcommand{\EM}[1]{{\it\textcolor{Maroon}{#1}}}
\newcommand{\BF}[1]{{\boldmath{\bf #1}\unboldmath}}
\newcommand{\B}{\{0,1\}}
\newcommand{\ONE}{\mathbf{1}}
\newcommand{\ZERO}{\mathbf{0}}
\newcommand{\G}{\mathbb{G}}
\newcommand{\cst}{\mathrm{cst}}
\newcommand{\id}{\mathrm{id}}

\title{Interaction graphs of isomorphic automata networks I:\\ complete digraph and minimum in-degree}

\date{\today}

\author{
Florian Bridoux\footnote{\scriptsize Universit\'e Côte d’Azur, CNRS, I3S, Sophia Antipolis, France. ({\tt bridoux@i3s.unice.fr})},
Kévin Perrot\footnote{\scriptsize Aix-Marseille Universit\'e, Universit\'e de Toulon, CNRS, LIS, Marseille, France. ({\tt kevin.perrot@lis-lab.fr})},
Aymeric Picard Marchetto\footnote{\scriptsize Universit\'e Côte d’Azur, CNRS, I3S, Sophia Antipolis, France. ({\tt  picard@i3s.unice.fr})},
Adrien Richard\footnote{\scriptsize Universit\'e Côte d’Azur, CNRS, I3S, Sophia Antipolis, France. ({\tt  adrien.richard@cnrs.fr})},
}

\begin{document}

\maketitle

\begin{abstract}
  An automata network with $n$ components over a finite alphabet $Q$ of size $q$
  is a discrete dynamical system described by the successive iterations of a
  function $f:Q^n\to Q^n$. In most applications, the main parameter is the
  interaction graph of $f$: the digraph with vertex set $[n]$ that contains an
  arc from $j$ to $i$ if $f_i$ depends on input $j$. What can be said on the
  set $\G(f)$ of the interaction graphs of the automata networks isomorphic to
  $f$? It seems that this simple question has never been studied. Here, we
  report some basic facts. First, we prove that if $n\geq 5$ or $q\geq 3$ and
  $f$ is neither the identity nor constant, then $\G(f)$ always contains the
  complete digraph $K_n$, with $n^2$ arcs. Then, we prove that $\G(f)$ always
  contains a digraph whose minimum in-degree is bounded as a function of $q$.
  Hence, if $n$ is large with respect to $q$,  then $\G(f)$ cannot only contain
  $K_n$. However, we prove that $\G(f)$ can contain only dense digraphs, with
  at least $\lfloor n^2/4 \rfloor$ arcs. 
\end{abstract}

\section{Introduction}

An \EM{automata network} with $n$ components over a finite alphabet $Q$ is a function 
\[
f:Q^n\to Q^n,\quad x=(x_1,\dots,x_n)\mapsto f(x)=(f_1(x),\dots,f_n(x)).
\]
The components $f_i:Q^n\to Q$, are usually called the \EM{local transition
functions} of the network, while $f$ is referred as the \EM{global transition
function} (but here we identify this function as the network). Automata
networks are also called \EM{finite dynamical systems}, and in the binary case,
when $Q=\{0,1\}$, the term \EM{Boolean networks} is applied. 

\medskip
The dynamics described by $f$ is explicitly represented by the digraph
\EM{$\Gamma(f)$} with vertex set $Q^n$ and an arc from $x$ to $f(x)$ for every
$x\in Q^n$. Two automata networks $f,h:Q^n\to Q^n$ are \EM{isomorphic}, $f\sim
h$ in notation, if $\Gamma(f)$ and $\Gamma(h)$ are isomorphic in the usual
sense; equivalently, there is permutation $\pi$ of $Q^n$ such that $\pi \circ f
= h \circ \pi$. An equivalence class of $\sim$ then corresponds to an unlabeled
digraph with $|Q|^n$ vertices in which each vertex has out-degree exactly one. 

\medskip
Automata networks have many applications. In particular, they are omnipresent
in the modeling of neural and gene networks (see \cite{N15} for a review). The
``network'' terminology comes from the fact that the \EM{interaction graph} of
$f$ is often considered as the main parameter: it is the digraph \EM{$G(f)$}
with vertex set $[n]=\{1,\dots,n\}$ and such that, for all $i,j\in [n]$, there
is an arc from $j$ to $i$ if $f_i$ depends on input $j$, that is, if there
exist $x,y\in Q^n$ which only differ in $x_j\neq y_j$ such that $f_i(x)\neq
f_i(y)$. 

\begin{figure}
  \centering
  \[
    \begin{array}{l}
      \begin{array}{cccclcc}
        \begin{array}{c|c}
          x & f(x)\\\hline
          000&000\\
          001&100\\
          010&000\\
          011&100\\
          100&000\\
          101&100\\
          110&000\\
          111&100
        \end{array}
        &&
        \begin{array}{c}
          \Gamma(f)\\
          \begin{tikzpicture}
            \useasboundingbox (-1.5,-0.6) rectangle (1.5,3);
            \node (000) at (0,0){$000$};
            \node (100) at (0,1.2){$100$};
            \node (010) at (-1.2,1.2){$010$};
            \node (110) at (1.2,1.2){$110$};
            \node (001) at (-1.2,2.4){$001$};
            \node (101) at (-0.4,2.4){$101$};
            \node (011) at (0.4,2.4){$011$};
            \node (111) at (1.2,2.4){$111$};
            \draw[->,thick] (000.-112) .. controls (-0.5,-0.7) and (0.5,-0.7) .. (000.-68);
            \path[thick,->,draw,black]
            (100) edge (000)
            (010) edge (000)
            (110) edge (000)
            (001) edge (100)
            (101) edge (100)
            (011) edge (100)
            (111) edge (100)
            ;
          \end{tikzpicture}
        \end{array}
        &&
        \begin{array}{l}
          f_1(x)=x_3\\[2mm]
          f_2(x)=0\\[2mm]
          f_3(x)=0\\
        \end{array}
        &&
        \begin{array}{c}
          G(f)\\[2mm]
          \begin{tikzpicture}
            \useasboundingbox (-0.7,-1.4) rectangle (2.2,0.4);
            \node[outer sep=1,inner sep=2,circle,draw,thick] (1) at (0,0){$1$};
            \node[outer sep=1,inner sep=2,circle,draw,thick] (2) at (1.5,0){$2$};
            \node[outer sep=1,inner sep=2,circle,draw,thick] (3) at (0.75,-1){$3$};
            \path[->,thick]
            (3) edge (1)
            ;
          \end{tikzpicture}
        \end{array}
        \\~\\
        \begin{array}{c|c}
          x & h(x)\\\hline
          000&000\\
          001&111\\
          010&111\\
          011&111\\
          100&111\\
          101&000\\
          110&000\\
          111&000
        \end{array}
        &&
        \begin{array}{c}
          \Gamma(h)\\
          \begin{tikzpicture}
            \useasboundingbox (-1.5,-0.6) rectangle (1.5,3);
            \node (000) at (0,0){$000$};
            \node (001) at (0,1.2){$111$};
            \node (010) at (-1.2,1.2){$101$};
            \node (011) at (1.2,1.2){$110$};
            \node (100) at (-1.2,2.4){$100$};
            \node (101) at (-0.4,2.4){$010$};
            \node (110) at (0.4,2.4){$001$};
            \node (111) at (1.2,2.4){$011$};
            \draw[->,thick] (000.-112) .. controls (-0.5,-0.7) and (0.5,-0.7) .. (000.-68);
            \path[thick,->,draw,black]
            (001) edge (000)
            (010) edge (000)
            (011) edge (000)
            (100) edge (001)
            (101) edge (001)
            (110) edge (001)
            (111) edge (001)
            ;
          \end{tikzpicture}
        \end{array}
        &&
        \begin{array}{l}
          h_1(x)=x_1 + (x_2 \lor x_3)\\[2mm]
          h_2(x)=x_1 + (x_2 \lor x_3)\\[2mm]
          h_3(x)=x_1 + (x_2 \lor x_3)\\
        \end{array}
        &&
        \begin{array}{c}
          G(h)\\[2mm]
          \begin{tikzpicture}
            \useasboundingbox (-0.7,-1.4) rectangle (2.2,0.4);
            \node[outer sep=1,inner sep=2,circle,draw,thick] (1) at (0,0){$1$};
            \node[outer sep=1,inner sep=2,circle,draw,thick] (2) at (1.5,0){$2$};
            \node[outer sep=1,inner sep=2,circle,draw,thick] (3) at (0.75,-1){$3$};
            \draw[->,thick] (1.{180+20}) .. controls (-0.8,-0.7) and (-0.8,+0.7) .. (1.{180-20});
            \draw[->,thick] (2.{0+20}) .. controls ({1.5+0.8},+0.7) and ({1.5+0.8},-0.7) .. (2.{0-20});
            \draw[->,thick] (3.{270+20}) .. controls ({0.75+0.7},{-1-0.8}) and ({0.75-0.7},{-1-0.8}) .. (3.{270-20});
            \path[->,thick]
            (2) edge[bend right=15] (1)
            (1) edge[bend right=15] (2)
            (3) edge[bend right=15] (1)
            (1) edge[bend right=15] (3)
            (2) edge[bend right=15] (3)
            (3) edge[bend right=15] (2)
            ;
          \end{tikzpicture}
        \end{array}
      \end{array}
    \end{array}
    \]
    \caption{\label{fig:isoexample}
    Two isomorphic automata networks on $\B^3$ with different interaction graphs.}
  \end{figure}
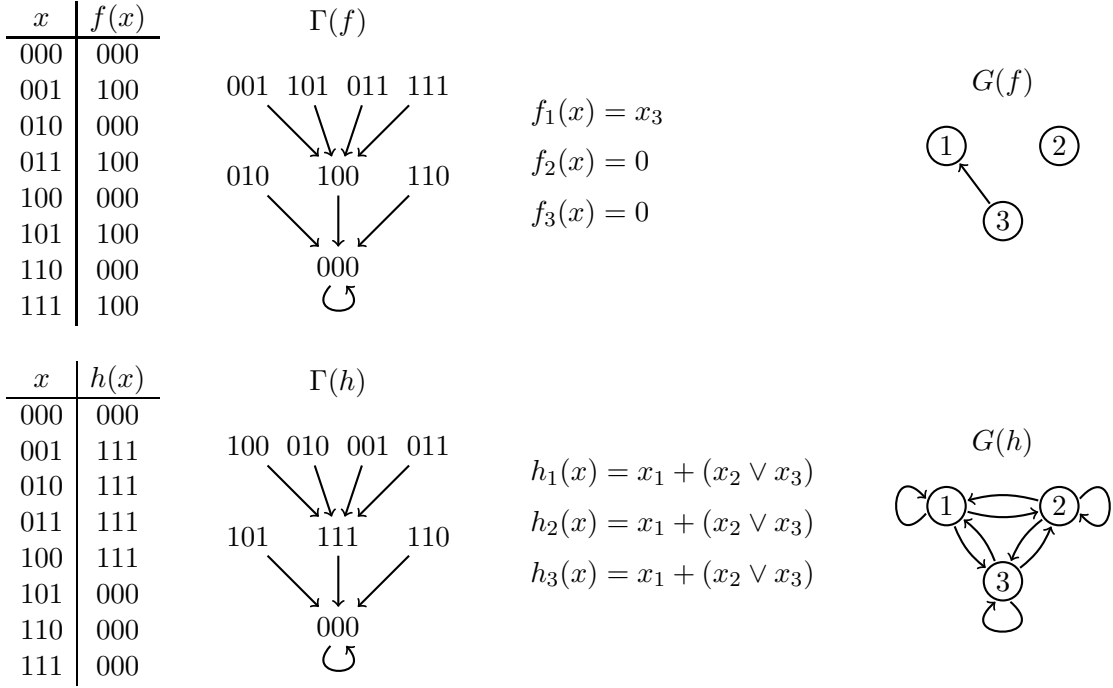

\medskip
In many applications, as modeling of gene networks, the interaction graph is
often well approximated while the actual dynamics is not \cite{TK01,N15}. One
is thus faced with the following question: {\em what can be said on the
dynamics described by $f$ from $G(f)$ only?} There are many results in this
direction; see \cite{G20} for a review. In most cases, the studied dynamical
properties are invariant by isomorphism: number of fixed points, number of
images, number of periodic configurations, number and size of limit cycles,
transient length and so on. However, the interaction graph is {\em not}
invariant by isomorphism: even if $f$ and $h$ are isomorphic, their interaction
graphs can be very different, as shown by Figure \ref{fig:isoexample}. This
variaton can give some {\em limit} to the central question stated above and we
think it deserves some study.

\medskip
Hence, we propose a systematic study of this variation, considering the set
\EM{$\G(f)$} of interaction graphs of automata networks isomorphic to $f$:
\[
\G(f)=\{G(h)\mid h:Q^n\to Q^n,~h\sim f\}.
\]
If $G\in\G(f)$ we say that $f$ \EM{can be produced} by $G$. For instance, if
$f$ is constant, we write this $f=\cst$, then $\G(f)$ contains a single
digraph, the digraph on $[n]$ without arcs, and if $f$ is the identity, we
write this $f=\id$, then $\G(f)$ also contains a single digraph, the digraph on
$[n]$ with $n$ loops (cycles of length one) and no other arcs. So the identity
and constant automata networks can be produced by a single digraph.

\medskip
Even if the systematic study of $\G(f)$ seems to us a new perspective, some
classical results can be stated in this setting. For instance, let us say that
a digraph is a \EM{companion digraph} if it can be obtained from a directed
path with $n$ vertices by adding between $0$ and $n$ arcs, all pointing to the
first vertex of the path. If $f$ is a \EM{linear network} with $n$ components
over the field  $Q=\mathbb{F}_p$ for some prime $p$, that is if $f(x)=Mx$ for
some $n\times n$ matrices $M$ over $\mathbb{F}_p$, then considering the
Frobenius normal form of $M$, we obtain that $f$ is isomorphic to a linear
network $h$ whose interaction graph $G$ is a disjoint union of companion
digraphs \cite{E59}. Consequently, any linear network $f$ can be produced by a
digraph $G\in\G(f)$ which is a disjoint union of companion digraphs, and which
is thus very sparse: $G$ has at most $2n-1$ arcs.

\medskip
The results we give on $\G(f)$ in this paper differ in that they do not rely on
any hypothesis made on $f$: we work with the whole set \EM{$F(n,q)$} of
automata networks with $n$ components over an alphabet of size $q$. Our first
result shows that, excepted few exceptions, $f$ can always be produced by the
\EM{complete digraph} on $[n]$, that is the digraph on $[n]$ with $n^2$ arcs,
denoted \EM{$K_n$}.

\begin{theorem}\label{thm:complete}
  Let $n,q\geq 2$ with  $n\geq 5$ or $q\geq 3$. If $f\in F(n,q)$ and
  $f\neq\cst,\id$, then~$K_n\in\G(f)$.
\end{theorem}

This indeed gives some limitations of what can be deduced on $f$ from $G(f)$:
if $G(f)=K_n$ {\em nothing} can be said on the dynamics of $f$ up to
isomorphism except that $f$ has at least two images and at most $q^n-1$ fixed
points. As explained above, if $f=\cst$ or $f=\id$ then $f$ can only be
produced by one digraph: the empty digraph and the digraph with $n$ loops,
respectively. The previous theorem says that if $f\neq\cst,\id$ can only be
produced by one digraph, then this digraph is necessarily $K_n$. But is there
at least one such automata network $f$, whose dynamics is completely specific
to the complete digraph? For $n=2$ and $q\geq 3$ we exhibit such $f$, and we
show that it does not exist when $n$ is large enough with respect to $q$. 

\medskip
We obtain these results by studying the minimum in-degree of the digraphs in
$\G(f)$. We first prove, using a difficult result in additive number theory,
that $\G(f)$ always contains a digraph whose minimum in-degree is bounded by a
constant $c_q$ that only depends on $q$, and which is $5$ for $q=2$. Hence, any
automata network $f$ is isomorphic to an automata network $h$ which contains a
``simple'' local transition function $h_i$, in that $h_i$ only depends on at most $c_q$
components. Conversely, we prove that, for some $f\in F(n,q)$, every digraph in
$\G(f)$ has minimum in-degree at least $2$. To state these results compactly,
for $n,q\geq 2$ we denote by \EM{$\delta^-(n,q)$} the minimum number $d$ such
that, for any $f\in F(n,q)$, some digraph in $\G(f)$ has minimum in-degree at
most $d$. 

\begin{theorem}\label{thm:in-degree}
  For all $n,q\geq 2$ with $n\geq 3$ or $q\geq 3$ we have 
  \[
    2\leq \delta^-(n,2)\leq 5
    \qquad\textrm{and}\qquad
    2\leq \delta^-(n,q)\leq (2+o(1))q
  \]
  where $o(1)$ tends to $0$ as $q$ tends to infinity,
  and does not depend on $n$.
\end{theorem}

Hence, for any $q\geq 3$, since $\delta^-(2,q)\geq 2$ there is $f\in F(2,q)$
such that every digraph $\G(f)$ has minimum in-degree at least $2$, and thus
$\G(f)=\{K_2\}$. So there are indeed some automata networks in $F(2,q)$ which
are only produced by $K_2$. On the other hand, for $n$ large enough with
respect to $q$, we have $\delta^-(n,q)<n$, so for any $f\in F(n,q)$ there is a
digraph $G\in\G(f)$ with minimum in-degree at most $n-1$, which is thus
distinct from $K_n$: no automata network in $F(n,q)$ is only produced by $K_n$.
That some automata networks in $F(n,q)$ are only produced by $K_n$ when $q$ is
large enough with respect to $n$ is open.  

\medskip
Even if for $n$ large with respect to $q$ there are no automata networks only
produced by $K_n$, we may ask if there are automata networks which are only
produced by dense digraphs. We answer this positively (and thus the fact,
mentioned above, that a linear network can always be produced by a sparse
digraph is specific to the linear case). 

\begin{theorem}\label{thm:dense}
  For all $n,q\geq 2$, there exists $f\in F(n,q)$ such that every digraph in
  $\G(f)$ has at least $\lfloor n^2/4 \rfloor$ arcs.
\end{theorem}

The previous results shed light on the case $|\G(f)|=1$. We may take the
converse direction and focus on automata networks $f$ for which $|\G(f)|$ is
large. Certainly $|\G(f)|\leq 2^{n^2}$ and the inequality is strict, since
otherwise $\G(f)$ contains the empty digraph, so that $f=\cst$ and thus
$|\G(f)|=1$. So $|\G(f)|\leq 2^{n^2}-1$. In a subsequent paper we will prove,
in particular, that this bound can be reached. Hence, while $K_n$ is in some
sense universal in that any automata network $f\neq\cst,\id$ with $n$
components can be produced by $K_n$, there are also universal automata networks
$f$ with $n$ components in that $f$ can be produced by all the digraphs on
$[n]$ except the empty one.

\medskip
The paper is organized as follows. Theorems \ref{thm:complete},
\ref{thm:in-degree} and \ref{thm:dense} are proved in Sections
\ref{sec:complete}, \ref{sec:in-degree} and \ref{sec:dense} respectively.
Concluding remarks are given in Section \ref{sec:concluding}.

\paragraph{Notations and terminologies} 

$n$ and $q$ are always integers greater than $1$. We set
$\EM{[n]}=\{1,\dots,n\}$ and $\EM{\Q{q}} = \{0,1,\dots,q-1\}$. Elements of
$\Q{q}^n$ are called \EM{configurations}, and those of $[n]$ are called
\EM{components}. We denote by \EM{$F(n,q)$} the set of functions
$f:\Q{q}^n\to\Q{q}^n$; so we always assume, without loss, that $Q=\Q{q}$. In
the following, the members of $F(n,q)$ are referred as functions instead of
automata networks. We denote by \EM{$\ZERO$} (resp. \EM{$\ONE$}) the
configurations $x$ such that $x_i = 0$ (resp. $1$) for all $i \in [n]$. For
$k\in\Q{q}$, we denote \EM{$ke_i$} the configuration $x$ such that $x_i=k$ and
$x_j=0$ for all $j\neq i$. We write \EM{$e_i$} instead of \EM{$1e_i$}. The
\EM{(Hamming) weight} of $x\in\Q{q}^n$ is the number of $i\in [n]$ with
$x_i\neq 0$. For $x,y \in \Q{q}^n$, the sum $\EM{x+y}$ is applied
component-wise modulo $q$. Hence, $x$ and $x+e_i$ differ only in component $i$.
So the interaction graph $G(f)$ of $f\in F(n,q)$ has an arc from $j$ to $i$ if
and only if there is $x\in\Q{q}^n$ such that $f_i(x)\neq f_i(x+e_j)$. A
\EM{fixed point} is a configuration $x$ such that $f(x)=x$. A \EM{limit cycle}
of $f$ is a cycle of $\Gamma(f)$. Hence, fixed points correspond to limit
cycles of length one. An \EM{independent set} of $f$ is an independent set of
$\Gamma(f)$, equivalently, it is a set $A\subseteq \Q{q}^n$ such that $f(A)\cap
A=\emptyset$. 

\section{Complete digraph}\label{sec:complete}

In this section, we will prove Theorem \ref{thm:complete} which states that
$K_n \in \G(f)$ for every $f \in F(n,q)$ with $f\neq\cst,\id$ and $n\geq 5$ or
$q\geq 3$. For that, we exhibit several sufficient conditions for
$K_n\in\G(f)$, and we then prove that at least one of these sufficient
conditions must be satisfied. As the proof for $q \geq 3$ is much simpler, we
will start with it. We will then treat the case $n\geq 5$ independently. 

\subsection{\BF{Case $q\geq 3$}}

We say that $f\in F(n,q)$ contains a \EM{complete-pattern} if there is $x\in
\Q{q}^n$ and a subset $A\subseteq \Q{q}^n$ of size $n$ such that $x,f(x)\notin
A\cup f(A)$ and either $x\neq f(x)$ or $A\cap f(A)=\emptyset$; in that case we
say that the complete-pattern is \EM{rooted} in $(x,A)$. 

\medskip
We prove, in Lemma~\ref{lem:complete-pattern_Kn} below, that, for $q\geq 3$ and
$n\geq 2$, we have $K_n\in \G(f)$ whenever $f$ has a complete-pattern. We then
prove, in Lemma \ref{lem:complete_pattern_exists}, that if $f\neq\cst,\id$ then
$f$ must contain a complete-pattern. Together, this proves Theorem
\ref{thm:complete} for $q \geq 3$.

\begin{lemma}\label{lem:complete-pattern_Kn}
  Let $f\in F(n,q)$ with $q \geq 3$ and $n\geq 2$. If $f$ has a
  complete-pattern then $K_n\in\G(f)$. 
\end{lemma}

\begin{figure}[h]
  \centering
  \begin{tikzpicture}
    \node (00) at (-1.0,1.0){$00$};
    \node (10) at (0,1.0){$10$};
    \node (01) at (1.0,1.0){$01$};
    \draw [draw=black,dotted] (-0.25,0.8) rectangle (1.25,1.2);
    \node (x) at (-1.0,1.4){\scriptsize$x$};
    \node (A) at (0.5,1.4){\scriptsize$A$};
    \node (11) at (-1.0,0){$11$};
    \node (20) at (0,0){$20$};
    \node (02) at (1.0,0){$02$};
    \path[thick,->,draw,black]
    (00) edge (11)
    (10) edge (20)
    (01) edge (02)
    ;
  \end{tikzpicture}
  \qquad
  \begin{tikzpicture}
    \node (00) at (-1.0,1.0){$00$};
    \node (10) at (0,1.0){$10$};
    \node (01) at (1.0,1.0){$01$};
    \draw [draw=black,dotted] (-0.25,0.8) rectangle (1.25,1.2);
    \node (x) at (-1.0,1.4){\scriptsize$x$};
    \node (A) at (0.5,1.4){\scriptsize$A$};
    \node (11) at (-1.0,0){$11$};
    \node (20) at (0.6,0){$20$};
    \path[thick,->,draw,black]
    (00) edge (11)
    (10) edge (20)
    (01) edge (20)
    ;
  \end{tikzpicture}
  \caption{\label{fig:q3_case_2}
  Labelling by $\pi$ of complete-patterns in the second case of Lemma \ref{lem:complete-pattern_Kn} for $n=2$.}
\end{figure}
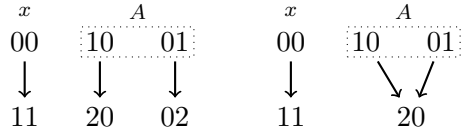

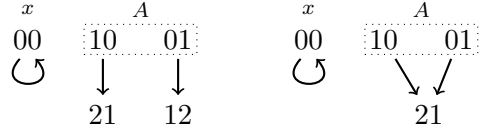
\begin{figure}[h]
  \centering
  \begin{tikzpicture}
    \node (00) at (-1.0,1.0){$00$};
    \node (10) at (0,1.0){$10$};
    \node (01) at (1.0,1.0){$01$};
    \draw [draw=black,dotted] (-0.25,0.8) rectangle (1.25,1.2);
    \node (x) at (-1.0,1.4){\scriptsize$x$};
    \node (A) at (0.5,1.4){\scriptsize$A$};
    \node (21) at (0,0){$21$};
    \node (12) at (1.0,0){$12$};
    \draw[->,thick] (00.-112) .. controls ({-1.0-0.5},{1.0-0.7}) and ({-1.0+0.5},{1.0-0.7}) .. (00.-68);
    \path[thick,->,draw,black]
    (10) edge (21)
    (01) edge (12)
    ;
  \end{tikzpicture}
  \qquad
  \begin{tikzpicture}
    \node (00) at (-1.0,1.0){$00$};
    \node (10) at (0,1.0){$10$};
    \node (01) at (1.0,1.0){$01$};
    \draw [draw=black,dotted] (-0.25,0.8) rectangle (1.25,1.2);
    \node (x) at (-1.0,1.4){\scriptsize$x$};
    \node (A) at (0.5,1.4){\scriptsize$A$};
    \node (21) at (0.6,0){$21$};
    \draw[->,thick] (00.-112) .. controls ({-1.0-0.5},{1.0-0.7}) and ({-1.0+0.5},{1.0-0.7}) .. (00.-68);
    \path[thick,->,draw,black]
    (10) edge (21)
    (01) edge (21)
    ;
  \end{tikzpicture}
  \caption{\label{fig:q3_case_3}
  Labelling by $\pi$ of complete-patterns in the third case of Lemma \ref{lem:complete-pattern_Kn} for $n=2$.}
\end{figure}

\begin{proof}
  Suppose that $f$ has a complete-pattern rooted in $(x,A)$. We have three cases:
  \begin{itemize}
    \item 
      Case 1: $x\neq f(x)$ and $f(A)\cap A\neq\emptyset$. Let $\pi$ be any permutation of $\Q{q}^n$ such that:
      \begin{itemize}
        \item
          $\pi(x)=\ZERO$,
        \item
          $\pi(f(x))=\ONE$,
        \item
          $\pi(A)=\{2e_1,\dots,2e_n\}$,
        \item
          $\pi(f(A)\setminus A) \subseteq \{0,2\}^n\setminus\{\ZERO,2e_1,\dots,2e_n\}$.
      \end{itemize}
      See Figure \ref{fig:q3_case_1} for an illustration. Note that since
      $f(A)\cap A\neq\emptyset$ the size of $f(A)\setminus A$ is at most $n-1$,
      and since $n-1\leq 2^n-(n+1)$ the configurations in $f(A)\setminus A$ can
      indeed have distinct images by $\pi$ inside
      $\{0,2\}^n\setminus\{\ZERO,2e_1,\dots,2e_n\}$. Let $h=\pi\circ f\circ
      \pi^{-1}$ and let us prove that $G(h)=K_n$. For all $i\in [n]$, the
      configurations $2e_i$ and $\ZERO$ differ only in component $i$, while
      $h(2e_i)\in\{0,2\}^n$ and $h(\ZERO) = \ONE$ differ in every component.
      Hence, $G(h)$ has an arc from $i$ to every component. Thus, $G(h) = K_n$
      and since $h\sim f$ we have $K_n\in\G(f)$.
    \item 
      Case 2: $x\neq f(x)$ and $f(A)\cap A=\emptyset$. Let $\pi$ be any permutation of $\Q{q}^n$ such that: 
      \begin{itemize}
        \item
          $\pi(x)=\ZERO$,
        \item
          $\pi(f(x))=\ONE$,
        \item
          $\pi(A)=\{e_1,\dots,e_n\}$,
        \item
          $\pi(f(A)) \subseteq \{2e_1,\dots,2e_n\}$.
      \end{itemize}
      See Figure \ref{fig:q3_case_2} for an illustration. Let $h=\pi\circ
      f\circ \pi^{-1}$ and let us prove that $G(h)=K_n$. For all $i\in [n]$, the
      configurations $e_i$ and $\ZERO$ differ only in component $i$, while
      $h(e_i)\in \{0,2\}^n$ and $h(\ZERO) = \ONE$ differ in every component.
      Hence, $G(h)$ has an arc from $i$ to every component. Thus, $G(h) = K_n$
      and since $h\sim f$ we have $K_n\in\G(f)$.
    \item 
      Case 3: $x= f(x)$ and $f(A)\cap A=\emptyset$. Let $\pi$ be any permutation of $\Q{q}^n$ such that:
      \begin{itemize}
        \item
          $\pi(x)=\ZERO$,
        \item
          $\pi(A)= \{e_1,\dots,e_n\}$,
        \item
          $\pi(f(A)) \subseteq \{\ONE+e_1,\dots,\ONE+e_n\}$.
      \end{itemize}
      See Figure \ref{fig:q3_case_3} for an illustration. Let $h=\pi\circ
      f\circ \pi^{-1}$ and let us prove that $G(h)=K_n$. For all $i\in [n]$, the
      configurations $e_i$ and $\ZERO$ differ only in component $i$, while
      $h(e_i)\in\{1,2\}^n$ and $h(\ZERO) = \ZERO$ differ in every component.
      Hence, $G(h)$ has an arc from $i$ to every component. Thus, $G(h) = K_n$
      and since $h\sim f$ we have $K_n\in\G(f)$.
  \end{itemize}
\end{proof}

\begin{figure}[h]
  \centering
  \begin{tikzpicture}
    \node (00) at (-1.0,1.0){$00$};
    \node (20) at (0,1.0){$20$};
    \node (02) at (1.0,1.0){$02$};
    \draw [draw=black,dotted] (-0.25,0.8) rectangle (1.25,1.2);
    \node (x) at (-1.0,1.4){\scriptsize$x$};
    \node (A) at (0.5,1.4){\scriptsize$A$};
    \node (11) at (-1.0,0){$11$};
    \draw[->,thick] (20.-112) .. controls ({0-0.5},{1.0-0.7}) and ({0+0.5},{1.0-0.7}) .. (20.-68);
    \draw[->,thick] (02.-112) .. controls ({1.0-0.5},{1.0-0.7}) and ({1.0+0.5},{1.0-0.7}) .. (02.-68);
    \path[thick,->,draw,black]
    (00) edge (11)
    ;
  \end{tikzpicture}
  \quad
  \begin{tikzpicture}
    \node (00) at (-1.0,1.0){$00$};
    \node (20) at (0,1.0){$20$};
    \node (02) at (1.0,1.0){$02$};
    \draw [draw=black,dotted] (-0.25,0.8) rectangle (1.25,1.2);
    \node (x) at (-1.0,1.4){\scriptsize$x$};
    \node (A) at (0.5,1.4){\scriptsize$A$};
    \node (11) at (-1.0,0){$11$};
    \draw[->,thick] (02.-112) .. controls ({1.0-0.5},{1.0-0.7}) and ({1.0+0.5},{1.0-0.7}) .. (02.-68);
    \path[thick,->,draw,black]
    (00) edge (11)
    (20) edge (02)
    ;
  \end{tikzpicture}
  \quad
  \begin{tikzpicture}
    \node (00) at (-1.0,1.0){$00$};
    \node (20) at (0,1.0){$20$};
    \node (02) at (1.0,1.0){$02$};
    \draw [draw=black,dotted] (-0.25,0.8) rectangle (1.25,1.2);
    \node (x) at (-1.0,1.4){\scriptsize$x$};
    \node (A) at (0.5,1.4){\scriptsize$A$};
    \node (11) at (-1.0,0){$11$};
    \node (22) at (0,0){$22$};
    \draw[->,thick] (02.-112) .. controls ({1.0-0.5},{1.0-0.7}) and ({1.0+0.5},{1.0-0.7}) .. (02.-68);
    \path[thick,->,draw,black]
    (00) edge (11)
    (20) edge (22)
    ;
  \end{tikzpicture}
  \quad
  \begin{tikzpicture}
    \node (00) at (-1.0,1.0){$00$};
    \node (20) at (0,1.0){$20$};
    \node (02) at (1.0,1.0){$02$};
    \draw [draw=black,dotted] (-0.25,0.8) rectangle (1.25,1.2);
    \node (x) at (-1.0,1.4){\scriptsize$x$};
    \node (A) at (0.5,1.4){\scriptsize$A$};
    \node (11) at (-1.0,0){$11$};
    \path[thick,->,draw,black]
    (00) edge (11)
    (20) edge[bend left=15] (02)
    (02) edge[bend left=15] (20)
    ;
  \end{tikzpicture}
  \quad
  \begin{tikzpicture}
    \node (00) at (-1.0,1.0){$00$};
    \node (20) at (0,1.0){$20$};
    \node (02) at (1.0,1.0){$02$};
    \draw [draw=black,dotted] (-0.25,0.8) rectangle (1.25,1.2);
    \node (x) at (-1.0,1.4){\scriptsize$x$};
    \node (A) at (0.5,1.4){\scriptsize$A$};
    \node (11) at (-1.0,0){$11$};
    \node (22) at (0,0){$22$};
    \path[thick,->,draw,black]
    (00) edge (11)
    (20) edge (22)
    (02) edge (20)
    ;
  \end{tikzpicture}
  \caption{\label{fig:q3_case_1}
  Labelling by $\pi$ of complete-patterns in the first case of Lemma \ref{lem:complete-pattern_Kn} for $n=2$.
  }
\end{figure}
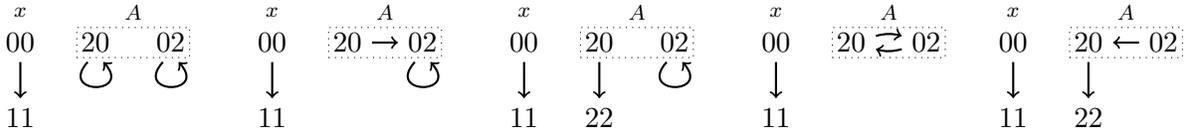

\begin{lemma}\label{lem:complete_pattern_exists}
  Let $f\in F(n,q)$ with $q \geq 3$ and $n\geq 2$. If $f\neq\cst,\id$ then $f$
  has a complete-pattern. 
\end{lemma}

\begin{proof}
  Suppose that there exists a configuration $a$ with at least $n+1$ pre-images
  by $f$. Then there is a set $A$ of size $n$ contained in
  $f^{-1}(a)\setminus\{a\}$. Since $f\neq\cst$, $f$ has an image $y\neq a$. Let
  $x\in f^{-1}(y)$. Then $x$ and $f(x)=y$ are not in in the union of $A$ and
  $f(A)=\{a\}$, and since $a\notin A$, $f$ has a complete-pattern rooted in
  $(x,A)$. 

  \medskip
  Suppose now that every configuration has at most $n$ pre-images by $f$. Since
  $f\neq\id$, there is $x\in\Q{q}^n$ such that $f(x)\neq x$. Let $B$ be the
  images of $f$ distinct from $x$ and $f(x)$. Then
  $|f^{-1}(B)|+|f^{-1}(x)|+|f^{-1}(f(x))|=q^n$ and we deduce that
  $|f^{-1}(B)|\geq q^n-2n$ which is at least $n+1$ since $q\geq 3$. So there is
  a subset $A\subseteq f^{-1}(B)\setminus\{f(x)\}$ of size $n$, and since
  $x,f(x)\notin A\cup f(A)$ we deduce that $f$ has a complete pattern rooted
  in $(x,A)$. 
\end{proof}

\subsection{\BF{Case $n\geq 5$}}

We proceed by showing that, given $f\in F(n,q)$ with $n\geq 5$, $q\geq 2$ and
$f\neq\cst,\id$, we have $K_n\in\G(f)$ if at least one of the following three
conditions holds: 
\begin{itemize}
  \item
    $f$ has at least $2n$ fixed points, 
  \item
    $f$ has at least $n$ limit cycles of length $\geq 3$,
  \item
    $f$ has an independent set of size $\geq 2n$. 
\end{itemize}
We then prove that, if $f\neq\cst,\id$, then at least one of the three
conditions must hold. Together this proves Theorem \ref{thm:complete} for
$n\geq 5$. We start by proving that the first two conditions, involving fixed
points and limit cycles, are sufficient for $K_n\in\G(f)$. 

\begin{lemma}\label{lem:fixed_points}
  Let $f\in F(n,q)$ with $f\neq\id$. If $f$ has at least $2n$ fixed points,
  then $K_n\in\G(f)$.  
\end{lemma}

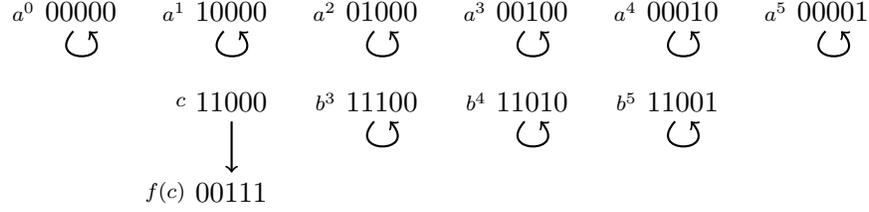
\begin{figure}[h]
  \centering
  \begin{tikzpicture}
    \node[label={[label distance=-4]180:{\scriptsize$a^0$}}] (00000) at (0,1.2){$00000$};
    \node[label={[label distance=-4]180:{\scriptsize$a^1$}}] (10000) at (2,1.2){$10000$};
    \node[label={[label distance=-4]180:{\scriptsize$a^2$}}] (01000) at (4,1.2){$01000$};
    \node[label={[label distance=-4]180:{\scriptsize$a^3$}}] (00100) at (6,1.2){$00100$};
    \node[label={[label distance=-4]180:{\scriptsize$a^4$}}] (00010) at (8,1.2){$00010$};
    \node[label={[label distance=-4]180:{\scriptsize$a^5$}}] (00001) at (10,1.2){$00001$};
    \node[label={[label distance=-4]180:{\scriptsize$c$}}] (11000) at (2,0){$11000$};
    \node[label={[label distance=-4]180:{\scriptsize$b^3$}}] (11100) at (4,0){$11100$};
    \node[label={[label distance=-4]180:{\scriptsize$b^4$}}] (11010) at (6,0){$11010$};
    \node[label={[label distance=-4]180:{\scriptsize$b^5$}}] (11001) at (8,0){$11001$};
    \node[label={[label distance=-4]180:{\scriptsize$f(c)$}}] (00111) at (2,-1.2){$00111$};
    \draw[->,thick] (00000.-112) .. controls (-0.5,{1.2-0.7}) and (0.5,{1.2-0.7}) .. (00000.-68);
    \draw[->,thick] (10000.-112) .. controls (1.5,{1.2-0.7}) and (2.5,{1.2-0.7}) .. (10000.-68);
    \draw[->,thick] (01000.-112) .. controls (3.5,{1.2-0.7}) and (4.5,{1.2-0.7}) .. (01000.-68);
    \draw[->,thick] (00100.-112) .. controls (5.5,{1.2-0.7}) and (6.5,{1.2-0.7}) .. (00100.-68);
    \draw[->,thick] (00010.-112) .. controls (7.5,{1.2-0.7}) and (8.5,{1.2-0.7}) .. (00010.-68);
    \draw[->,thick] (00001.-112) .. controls (9.5,{1.2-0.7}) and (10.5,{1.2-0.7}) .. (00001.-68);
    \draw[->,thick] (11100.-112) .. controls (3.5,-0.7) and (4.5,-0.7) .. (11100.-68);
    \draw[->,thick] (11010.-112) .. controls (5.5,-0.7) and (6.5,-0.7) .. (11010.-68);
    \draw[->,thick] (11001.-112) .. controls (7.5,-0.7) and (8.5,-0.7) .. (11001.-68);
    \path[thick,->,draw,black]
    (11000) edge (00111)
    ;
  \end{tikzpicture}
  \caption{\label{fig:fixed_points}
  Labelling by $\pi$ in Lemma \ref{lem:fixed_points} for $n = 5$.
  }

\end{figure}

\begin{proof}
  Since $f\neq \id$ we have $f(c)\neq c$ for some configuration $c$, and since
  $f$ has at least $2n$ fixed points, it has $2n-1$ fixed points distinct from
  $f(c)$, say $a^0,a^1,\dots,a^n,b^3,\dots,b^n$ (these are also distinct from
  $c$ since $c$ is not a fixed point). Let $\pi$ be any permutation of
  $\Q{q}^n$ such that:
  \begin{itemize}
    \item
      $\pi(a^0)=\ZERO$, 
    \item
      $\pi(a^i)=e_i$ for $1\leq i\leq n$, 
    \item
      $\pi(b^i)=e_1+e_2+e_i$ for $3\leq i\leq n$, 
    \item
      $\pi(c)=e_1+e_2$,
    \item
      $\pi(f(c))=e_1+e_2+\ONE$. 
  \end{itemize}
  See Figure \ref{fig:fixed_points} for an illustration. Let $h=\pi\circ f\circ
  \pi^{-1}$ and let us prove that $G(h)=K_n$. For $i\in [n]$, the
  configurations $\ZERO$ and $e_i$ only differ in component $i$, and since
  $h(\ZERO)=\ZERO$ and $h(e_i)=e_i$ we deduce that $G(h)$ has a loop on $i$. It
  remains to prove that $G(h)$ has an arc from $i$ to $j$ for distinct $i,j\in
  [n]$. The configurations $e_2$ and $e_1+e_2$ only differ in component $1$,
  while $h(e_2)=e_2$ and $h(e_1+e_2)=e_1+e_2+\ONE$ differ in every component
  $j\neq 1$. So $G(h)$ has an arc from $1$ to every component $j\neq 1$.
  Similarly, the configurations $e_1$ and $e_1+e_2$ only differ in component
  $2$, while $h(e_1)=e_1$ and $h(e_1+e_2)=e_1+e_2+\ONE$ differ in every
  component $j\neq 2$. So $G(h)$ has an arc from $2$ to every component $j\neq
  2$. Finally, for $3\leq i\leq n$, the configurations $e_1+e_2$ and
  $e_1+e_2+e_i$ only differ in component $i$, while
  $h(e_1+e_2+e_i)=e_1+e_2+e_i$ and $h(e_1+e_2)=e_1+e_2+\ONE$ differ in  every
  component $j\neq i$. So $G(h)$ has an arc from $i$ to every $j\neq i$.  Thus,
  $G(h) = K_n$ and since $h\sim f$ we have $K_n\in\G(f)$.
\end{proof}

\begin{lemma}\label{lem:limit_cycles}
  Let $f\in F(n,q)$. If $f$ has at least $n$ limit cycles of length $\geq 3$,
  then $K_n\in\G(f)$.  
\end{lemma}

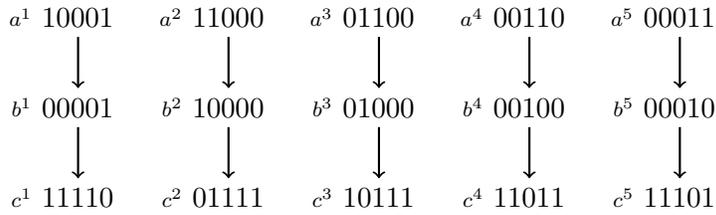
\begin{figure}[h]
  \centering
  \begin{tikzpicture}
    \useasboundingbox (-0.3,-0.6) rectangle (8.3,3);
    \node[label={[label distance=-4]180:{\scriptsize$a^1$}}] (10001) at (0,2.4){$10001$};
    \node[label={[label distance=-4]180:{\scriptsize$b^1$}}] (00001) at (0,1.2){$00001$};
    \node[label={[label distance=-4]180:{\scriptsize$c^1$}}] (11110) at (0,0){$11110$};
    \node[label={[label distance=-4]180:{\scriptsize$a^2$}}] (11000) at (2,2.4){$11000$};
    \node[label={[label distance=-4]180:{\scriptsize$b^2$}}] (10000) at (2,1.2){$10000$};
    \node[label={[label distance=-4]180:{\scriptsize$c^2$}}] (01111) at (2,0){$01111$};
    \node[label={[label distance=-4]180:{\scriptsize$a^3$}}] (01100) at (4,2.4){$01100$};
    \node[label={[label distance=-4]180:{\scriptsize$b^3$}}] (01000) at (4,1.2){$01000$};
    \node[label={[label distance=-4]180:{\scriptsize$c^3$}}] (10111) at (4,0){$10111$};
    \node[label={[label distance=-4]180:{\scriptsize$a^4$}}] (00110) at (6,2.4){$00110$};
    \node[label={[label distance=-4]180:{\scriptsize$b^4$}}] (00100) at (6,1.2){$00100$};
    \node[label={[label distance=-4]180:{\scriptsize$c^4$}}] (11011) at (6,0){$11011$};
    \node[label={[label distance=-4]180:{\scriptsize$a^5$}}] (00011) at (8,2.4){$00011$};
    \node[label={[label distance=-4]180:{\scriptsize$b^5$}}] (00010) at (8,1.2){$00010$};
    \node[label={[label distance=-4]180:{\scriptsize$c^5$}}] (11101) at (8,0){$11101$};
    \path[thick,->,draw,black]
    (10001) edge (00001)
    (00001) edge (11110)
    (11000) edge (10000)
    (10000) edge (01111)
    (01100) edge (01000)
    (01000) edge (10111)
    (00110) edge (00100)
    (00100) edge (11011)
    (00011) edge (00010)
    (00010) edge (11101)
    ;
  \end{tikzpicture}
  \caption{\label{fig:limit_cycles}
  Labelling by $\pi$ in Lemma \ref{lem:limit_cycles} for $n = 5$. 
  }
\end{figure}

\begin{proof}
  Suppose that $f$ has $n$ limit cycles of length $\geq 3$; this implies $n\geq
  4$. Then there are $n$ configurations $a^1,\dots,a^n$ that belong to distinct
  limit cycles of $f$ of length $\geq 3$. For $i\in [n]$, let $b^i=f(a^i)$ and
  $c^i=f(b^i)$. Then $a^1,\dots,a^n$, $b^1,\dots,b^n$, $c^1,\dots,c^n$ are all
  distinct. Let $\pi$ be any permutation of $\Q{q}^n$ such that 
  \begin{itemize}
    \item
      $\pi(a^i)=e_{i-1}+e_i$,
    \item
      $\pi(b^i)=e_{i-1}$,
    \item
      $\pi(c^i)=e_{i-1}+\ONE$,
  \end{itemize}
  where $e_0$ means $e_n$. See Figure \ref{fig:limit_cycles} for an
  illustration. Note that the map $\pi$ is indeed injective on the $3n$
  selected configurations: the weights of $\pi(a^i)$, $\pi(b^i)$, and
  $\pi(c^i)$ are $2$, $1$, and $n-1 \geq 4-1$ respectively. Let $h=\pi\circ
  f\circ \pi^{-1}$ and let us prove that $G(h)=K_n$. For $i\in [n]$, the
  configurations $e_{i-1}$ and $e_{i-1}+e_i$ only differ in component $i$,
  while $h(e_{i-1})=e_{i-1}+\ONE$ and  $h(e_{i-1}+e_i)=e_{i-1}$ differ in every
  component. So $G(h)$ has an arc from $i$ to every $j\in [n]$. Thus, $G(h) =
  K_n$ and since $h\sim f$ we have $K_n\in\G(f)$.
\end{proof}

We now prove that the third condition given at the beginning of the subsection,
involving independent sets, is sufficient for $K_n\in\G(f)$. The precise
statement is as follows. 

\begin{lemma}\label{lem:independent_set}
  Let $f\in F(n,q)$ with $n\geq 5$ and $f\neq\cst$. If $f$ has an independent
  set of size at least $2n$, then $K_n\in\G(f)$.  
\end{lemma}

The proof is more involved than for fixed points and limit cycles, and so we
split the proof. We first prove that, if $n\geq 5$ and $f\neq\cst$, then
$K_n\in\G(f)$ whenever at least one of the following condition holds:
\begin{itemize}
  \item $f$ has an independent set $A$ with $|A|\geq n+k$ and $|f(A)|=2k$ for
    some $1\leq k\leq n$.
  \item $f$ has an independent set $A$ with $|A|>n$ and $|f(A)|=1$.
\end{itemize}
We then prove that if $f$ has an independent set of size at least $2n$, then at
least one of the two conditions must hold, and this proves
Lemma~\ref{lem:independent_set}.

\begin{lemma}\label{lem:independent_set_2k}
  Let $f\in F(n,q)$ with $n\geq 5$. If $f$ has an independent set $A$ such that
  $|A|\geq n+k$ and $|f(A)|=2k$ for some $1\leq k\leq n$, then $K_n\in\G(f)$.  
\end{lemma}

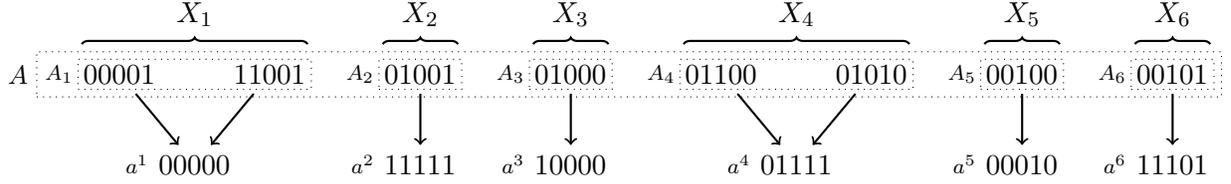
\begin{figure}[h]
  \centering
  \begin{tikzpicture}
    \def\c{2}
    \node (A11) at ({0*\c},1.2){$00001$};
    \node (A12) at ({1*\c},1.2){$11001$};
    \node (A2)  at ({2*\c},1.2){$01001$};
    \node (A3)  at ({3*\c},1.2){$01000$};
    \node (A41) at ({4*\c},1.2){$01100$};
    \node (A42) at ({5*\c},1.2){$01010$};
    \node (A5)  at ({6*\c},1.2){$00100$};
    \node (A6)  at ({7*\c},1.2){$00101$};

    \draw[decorate,decoration={brace},thick] ({(0*\c)-0.5},1.6) --  ({(1*\c)+0.5},1.6);
    \draw[decorate,decoration={brace},thick] ({(2*\c)-0.5},1.6) --  ({(2*\c)+0.5},1.6);
    \draw[decorate,decoration={brace},thick] ({(3*\c)-0.5},1.6) --  ({(3*\c)+0.5},1.6);
    \draw[decorate,decoration={brace},thick] ({(4*\c)-0.5},1.6) --  ({(5*\c)+0.5},1.6);
    \draw[decorate,decoration={brace},thick] ({(6*\c)-0.5},1.6) --  ({(6*\c)+0.5},1.6);
    \draw[decorate,decoration={brace},thick] ({(7*\c)-0.5},1.6) --  ({(7*\c)+0.5},1.6);

    \draw [draw=black,dotted] ({(0*\c)-1.1},{1.2+0.3}) rectangle  ({(7*\c)+0.65},{1.2-0.3});
    \node at ({(0*\c)-1.35},1.2){\small$A$};

    \draw [draw=black,dotted] ({(0*\c)-0.55},{1.2+0.2}) rectangle  ({(1*\c)+0.55},{1.2-0.2});
    \draw [draw=black,dotted] ({(2*\c)-0.55},{1.2+0.2}) rectangle  ({(2*\c)+0.55},{1.2-0.2});
    \draw [draw=black,dotted] ({(3*\c)-0.55},{1.2+0.2}) rectangle  ({(3*\c)+0.55},{1.2-0.2});
    \draw [draw=black,dotted] ({(4*\c)-0.55},{1.2+0.2}) rectangle  ({(5*\c)+0.55},{1.2-0.2});
    \draw [draw=black,dotted] ({(6*\c)-0.55},{1.2+0.2}) rectangle  ({(6*\c)+0.55},{1.2-0.2});
    \draw [draw=black,dotted] ({(7*\c)-0.55},{1.2+0.2}) rectangle  ({(7*\c)+0.55},{1.2-0.2});

    \node at ({(0*\c)-0.8},1.2){\scriptsize$A_1$};
    \node at ({(2*\c)-0.8},1.2){\scriptsize$A_2$};
    \node at ({(3*\c)-0.8},1.2){\scriptsize$A_3$};
    \node at ({(4*\c)-0.8},1.2){\scriptsize$A_4$};
    \node at ({(6*\c)-0.8},1.2){\scriptsize$A_5$};
    \node at ({(7*\c)-0.8},1.2){\scriptsize$A_6$};

    \node at ({\c/2},2){$X_1$};
    \node at ({2*\c},2){$X_2$};
    \node at ({3*\c},2){$X_3$};
    \node at ({4*\c+\c/2},2){$X_4$};
    \node at ({6*\c},2){$X_5$};
    \node at ({7*\c},2){$X_6$};

    \node[label={[label distance=-4]180:{\scriptsize$a^1$}}] (a1) at ({0*\c+\c/2},0){$00000$};
    \node[label={[label distance=-4]180:{\scriptsize$a^2$}}] (a2) at ({2*\c},0){$11111$};
    \node[label={[label distance=-4]180:{\scriptsize$a^3$}}] (a3) at ({3*\c},0){$10000$};
    \node[label={[label distance=-4]180:{\scriptsize$a^4$}}] (a4) at ({4*\c+\c/2},0){$01111$};
    \node[label={[label distance=-4]180:{\scriptsize$a^5$}}] (a5) at ({6*\c},0){$00010$};
    \node[label={[label distance=-4]180:{\scriptsize$a^6$}}] (a6) at ({7*\c},0){$11101$};

    \path[thick,->,draw,black]
    (A11) edge (a1)
    (A12) edge (a1)
    (A2)  edge (a2)
    (A3)  edge (a3)
    (A41) edge (a4)
    (A42) edge (a4)
    (A5)  edge (a5)
    (A6)  edge (a6)
    ;
  \end{tikzpicture}
  \caption{\label{fig:independent_set_2k}
  Labelling by $\pi$ in Lemma \ref{lem:independent_set_2k} for $n = 5$ and
  $k=3$. The independent set $A$ is of size $|A|=n+k=8$ and its image
  $f(A)=\{a^1,\dots,a^6\}$ is of size $2k=6$. Since the consecutive sizes of
  the set $A_1,\dots,A_6$ are $(2,1,1,2,1,1)$, we have to choose a partition
  $I$ of $[n]$ into $6$ parts $I_1,\dots,I_6$ whose consecutive sizes are
  $(1,1,0,2,0,1)$. We take
  $(I_1,I_2,I_3,I_4,I_5,I_6)=(\{1\},\{2\},\emptyset,\{3,4\},\emptyset,\{5\})$.
  We have to choose one member $j_2,j_4,j_6$ in $X_2,X_4,X_6$ respectively. We
  choose $j_2=2$, $j_4=3$ and $j_6=5$. For these choices, we obtain the set
  $X_1,\dots,X_6$ described in the figure, which are used to label the
  configurations in $A_1,\dots,A_6$. Finally, we have to choose a set $Y$ of
  $k=3$ configurations, such that $Y,Y+\ONE$ and $X=X_1\cup\cdots\cup X_6$ are
  pairwise disjoint. We choose $Y=\{00000,10000,00010\}$, which is used to
  label $a_1,a_3,a_5$, and $Y+\ONE$ is then used to label $a_2,a_4,a_6$
  accordingly. Note that: $11001\in X_1$ and $01001\in X_2$ only differ in
  component $1$; $01000\in X_1$ and $01001\in X_2$ only differ in component
  $2$; $01000\in X_3$ and $01100\in X_4$ only differ in component $3$;
  $01000\in X_3$ and $01010\in X_4$ only differ in component $4$; $00100\in
  X_5$ and $00101\in X_6$ only differ in component $5$.}
\end{figure}

\begin{proof}
  Suppose that $f$ has an independent set $A$ of $f$ such that $|A|\geq n+k$
  and $|f(A)|=2k$ for some $1\leq k\leq n$, and suppose that $|A|$ is minimal
  for those properties. 

  \begin{itemize}
    \item[(1)]
      {\em $|A|=n+k$.}

      Indeed, if $|A|>n+k$, then $|A|>2k$ and thus there are distinct $a,a'\in
      A$ with $f(a)=f(a')$. Hence, $A'=A\setminus \{a\}$ is an independent set
      of $f$ such that $|A'|=|A|-1\geq n+k$ and $|f(A')|=2k$, and this
      contradicts the minimality of $|A|$. This proves (1).  
  \end{itemize}

  Let us write $f(A)=\{a^1,\dots,a^{2k}\}$, and let $A_p=f^{-1}(a^p)\cap A$ for
  $p\in[2k]$. Hence, $\{A_1,\dots,A_{2k}\}$ is a partition of $A$. We now come
  to the technical part of the proof. 

  \begin{itemize}
    \item[(2)]
      {\em There are disjoint subsets $X_1,\dots,X_{2k}$ of $\B^n$ of size
      $|A_1|,\dots,|A_{2k}|$ respectively, such that, for all $i\in [n]$, there
      are $p\in [k]$ and $x\in X_{2p-1}$ with $x+e_i\in X_{2p}$.}

      Let $I_1,\dots,I_{2k}$ be a partition of $[n]$ (with some members
      possibly empty) such that, for $1\leq \ell\leq 2k$, the size of $I_\ell$
      is $|A_\ell|-1$ if $\ell$ is odd, and $|A_\ell|$ otherwise; such a
      partition exists since $|A_1|+\cdots+|A_{2k}|=n+k$ by (1). 

      For $p\in [k]$, select a configuration $x^{2p-1}\in\B^n$, a component
      $j_{2p}\in I_{2p}$, and let 
      \[
        \begin{array}{lll}
          X_{2p-1}&=&\{x^{2p-1}\}\cup\{x^{2p-1}+e_{j_{2p}}+e_i\mid i\in I_{2p-1}\},\\[1mm]
          X_{2p}&=&\{x^{2p-1}+e_i\mid i\in I_{2p}\}. 
        \end{array}
      \]
      Clearly, $|X_{2p-1}|=|I_{2p-1}|+1=|A_{2p-1}|$ and $|X_{2p}|=|I_{2p}|=|A_{2p}|$.

      Let $i\in [n]$ and let us prove that there are $p\in [k]$ and $x\in
      X_{2p-1}$ with $x+e_i\in X_{2p}$. We have $i\in I_\ell$ for some $1\leq
      \ell\leq 2k$. If $\ell=2p-1$ then, setting $x=x^{2p-1}+e_{j_{2p}}+e_i$,
      we have $x\in X_{2p-1}$ and $x+e_i=x^{2p-1}+e_{j_{2p}}\in X_{2p}$. If
      $\ell=2p$ then, setting $x=x^{2p-1}$, we have $x\in X_{2p-1}$ and
      $x+e_i\in X_{2p}$. Thus, fixing the components $j_{2p}$ arbitrarily,
      we only have to prove that we can choose the configurations $x^{2p-1}$
      so that the sets $X_1,\dots,X_{2k}$ are
      pairwise disjoint.

      This is obvious if $k=1$, since $X_1$ and $X_2$ are disjoint for any
      choice of $x^1$. If $k=2$, then, since $n\geq 5$, the sets
      $X_1,X_2,X_3,X_4$ are disjoint by taking $x^1=\ZERO$ and $x^3=\ONE$:
      configurations in $X_1$ have weight $0$ or $2$; configurations in $X_2$
      are of weight $1$; configurations in $X_3$ are of weight $n$ or $n-2$;
      and configurations in $X_4$ are of weight $n-1$.

      Suppose now that $k\geq 3$ and choose $x^{2p-1}=e_{j_{2p-2}}$ for all
      $p\in[k]$, where $j_0$ means~$j_{2k}$. Then each $X_{2p-1}$ only contains
      configurations of weight $1$ or $3$ and each $X_{2p}$ only contains
      configurations of weight $2$. Hence, given $1\leq p<q\leq k$, we have to
      prove that $X_{2p-1}\cap X_{2q-1}=\emptyset$ and $X_{2p}\cap
      X_{2q}=\emptyset$.

      Suppose that $x\in X_{2p-1}\cap X_{2q-1}$. If the weight of $x$ is $1$
      then we deduce that $x=e_{j_{2p-2}}=e_{j_{2q-2}}$ which is false since
      $p\neq q$. If the weight of $x$ is $3$ then $x_i=1$ for some $i\in
      I_{2p-1}$ while $y_i=0$ for all $y\in X_{2q-1}$, a contradiction. Thus,
      $X_{2p-1}\cap X_{2q-1}=\emptyset$.

      Suppose now that $x\in X_{2p}\cap X_{2q}$. Then
      $x=e_{j_{2p-2}}+e_{i_{2p}}=e_{j_{2q-2}}+e_{i_{2q}}$ for some $i_{2p}\in
      I_{2p}$ and $i_{2q}\in I_{2q}$. Thus, $i_{2p}=j_{2q-2}$ and
      $i_{2q}=j_{2p-2}$. Hence, $j_{2q-2}\in I_{2p}$, and since $p<q$ this
      implies $q=p+1$. Also, $j_{2p-2}\in I_{2q}$ and since $p<q$ this implies
      $p=1$ and $q=k$, but then $q\neq p+1$ since $k\geq 3$. Thus,
      $X_{2p}\cap X_{2q}=\emptyset$. Hence, the sets $X_1,\dots,X_{2k}$ are
      indeed pairwise disjoint. This proves (2). 
  \end{itemize}

  Let $X_1,\dots,X_{2k}$ be as in (2). Let $X=X_1\cup\dots\cup X_{2k}$, which
  is of size $n+k$. 

  \begin{itemize}
    \item[(3)]
      {\em There is $Y\subseteq\B^n$ with $|Y|\geq n$ such that $Y$, $Y+\ONE$
      and $X$ are pairwise disjoint.}

      Let $Z$ be any subset of $\B^n$ of size $2^{n-1}$ such that $Z$ and
      $Z+\ONE$ are disjoint. Let $Y$ be the set of $y\in Z$ such that
      $y,y+\ONE\notin X$. Clearly, $Y,Y+\ONE,X$ are pairwise disjoint, so it
      is sufficient to prove that $|Y|\geq n$. Setting $X^0=X\cap Z$ and
      $X^1=X\cap (Z+\ONE)$, we have $Y=Z\setminus (X^0\cup (X^1+\ONE))$ and
      thus $|Y|\geq |Z|-|X^0|-|X^1|$. Since $Z$ and $Z+\ONE$ are disjoint,
      $X^0$ and $X^1$ are disjoint subsets of $X$, so 
      \[
        |Y|\geq |Z|-|X|=2^{n-1}-(n+k)\geq 2^{n-1}-2n\geq n. 
      \] 
      where the last inequality holds since $n\geq 5$. 
  \end{itemize}

  Let $Y=\{y^1,\dots,y^k\}$ be a subset of $\B^n$ of size $k$ such that $Y$,
  $Y+\ONE$ and $X$ are pairwise disjoint, which exists by (3). Hence, there is
  a permutation $\pi$ of $\B^n$ such that, for all $p\in [k]$:
  \begin{itemize}
    \item $\pi(a^{2p-1})=y^p$,
    \item $\pi(a^{2p})=y^p+\ONE$,
    \item $\pi(A_{2p-1})=X_{2p-1}$,
    \item $\pi(A_{2p})=X_{2p}$.
  \end{itemize}
  See Figure \ref{fig:independent_set_2k} for an illustration. Let $h=\pi\circ
  f\circ\pi^{-1}$ and let us prove that $G(h)=K_n$. Since the sets
  $X_1,\dots,X_{2k}$ are as in (2), for every $i\in [n]$ there is $p\in [k]$
  and $x\in X_{2p-1}$ with $x+e_i\in X_{2p}$. Since $\pi^{-1}(x)\in A_{2p-1}$
  and $\pi^{-1}(x+e_i)\in A_{2p}$, we have $h(x)=\pi(a^{2p-1})=y^p$ and
  $h(x+e_i)=\pi(a^{2p})=y^p+\ONE$. Thus, $h(x)$ and $h(x+e_i)$ differ in every
  component, and so $G(h)$ has an arc from $i$ to every $j\in [n]$. Thus, $G(h)
  = K_n$ and since $h\sim f$ we have $K_n\in\G(f)$.
\end{proof}

\begin{lemma}\label{lem:independent_set_1}
  Let $f\in F(n,q)$ with $n\geq 5$ and $f\neq\cst$. If $f$ has an independent
  set $A$ such that $|A|>n$ and $|f(A)|=1$, then $K_n\in\G(f)$.  
\end{lemma}

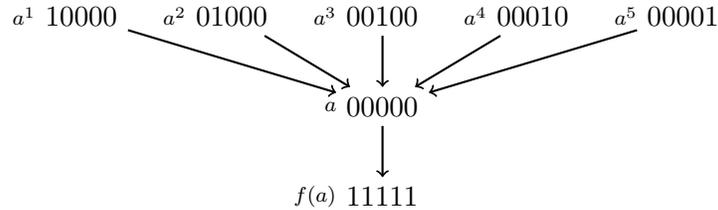
\begin{figure}[h]
  \centering
  \begin{tikzpicture}
    \def\c{2}
    \node[label={[label distance=-4]180:{\scriptsize$a$}}] (a) at ({2*\c},0){$00000$};
    \node[label={[label distance=-4]180:{\scriptsize$f(a)$}}] (fa) at ({2*\c},-1.2){$11111$};
    \node[label={[label distance=-4]180:{\scriptsize$a^1$}}] (a1) at ({0*\c},1.2){$10000$};
    \node[label={[label distance=-4]180:{\scriptsize$a^2$}}] (a2) at ({1*\c},1.2){$01000$};
    \node[label={[label distance=-4]180:{\scriptsize$a^3$}}] (a3) at ({2*\c},1.2){$00100$};
    \node[label={[label distance=-4]180:{\scriptsize$a^4$}}] (a4) at ({3*\c},1.2){$00010$};
    \node[label={[label distance=-4]180:{\scriptsize$a^5$}}] (a5) at ({4*\c},1.2){$00001$};
    \path[thick,->,draw,black]
    (a) edge (fa)
    (a1) edge (a)
    (a2) edge (a)
    (a3) edge (a)
    (a4) edge (a)
    (a5) edge (a)
    ;
  \end{tikzpicture}
  \caption{\label{fig:independent_set_1_case1}
  Labelling by $\pi$ in the first case of Lemma \ref{lem:independent_set_1} for $n=5$.
  }
\end{figure}

\begin{figure}[h]
  \centering
  \begin{tikzpicture}
    \def\c{2}
    \node[label={[label distance=-4]180:{\scriptsize$a$}}] (a) at ({2*\c},0){$11111$};
    \node[label={[label distance=-4]180:{\scriptsize$a^1$}}] (a1) at ({0*\c},1.2){$10000$};
    \node[label={[label distance=-4]180:{\scriptsize$a^2$}}] (a2) at ({1*\c},1.2){$01000$};
    \node[label={[label distance=-4]180:{\scriptsize$a^3$}}] (a3) at ({2*\c},1.2){$00100$};
    \node[label={[label distance=-4]180:{\scriptsize$a^4$}}] (a4) at ({3*\c},1.2){$00010$};
    \node[label={[label distance=-4]180:{\scriptsize$a^5$}}] (a5) at ({4*\c},1.2){$00001$};
    \node[label={[label distance=-4]180:{\scriptsize$b$}}] (b) at ({-1*\c},1.2){$00000$};
    \draw[->,thick] (b.-112) .. controls ({-1*\c-0.5},{1.2-0.7}) and ({-1*\c+0.5},{1.2-0.7}) .. (b.-68);
    \draw[->,thick] (a.-112) .. controls ({2*\c-0.5},{0-0.7}) and ({2*\c+0.5},{0-0.7}) .. (a.-68);
    \path[thick,->,draw,black]
    (a1) edge (a)
    (a2) edge (a)
    (a3) edge (a)
    (a4) edge (a)
    (a5) edge (a)
    ;
  \end{tikzpicture}
  \caption{\label{fig:independent_set_1_case2}
  Labelling by $\pi$ in the second case of Lemma \ref{lem:independent_set_1} for $n=5$.
  }
\end{figure}
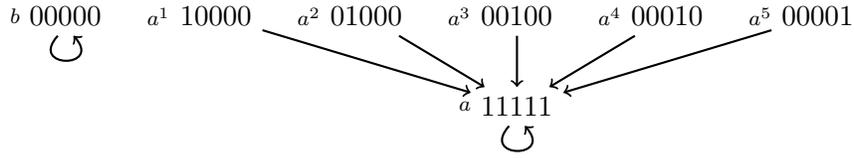

\begin{proof}
  Let $A$ be an independent set of $f$ with $|A|>n$ and $|f(A)|=1$. Then
  $f(A)=\{a\}$ for some configuration $a\in\Q{q}^n$.  Since $f\neq\cst$, there
  is $b\in\Q{q}^n$ with $f(b)\neq a$, and thus $b\notin A$.  We consider three
  cases. 
  \begin{itemize}
    \item
      Case 1: $f(a)\neq a$. Since $|A|>n$, there are $n$ configurations
      $a^1,\dots,a^n$ in $A$ distinct from $f(a)$. Then $a^1,\dots,a^n,a,f(a)$
      are all distinct, so there is a permutation $\pi$ of $\Q{q}^n$ such that:
      \begin{itemize}
        \item $\pi(a)=\ZERO$,
        \item $\pi(f(a))=\ONE$,
        \item $\pi(a^i)=e_i$ for $1\leq i\leq n$.
      \end{itemize}
      See Figure \ref{fig:independent_set_1_case1} for an illustration. Let
      $h=\pi\circ f\circ\pi^{-1}$ and let us prove that $G(h)=K_n$. For $i\in
      [n]$, $h(\ZERO)=\ONE$ and $h(e_i)=\ZERO$ differ in every component,
      therefore $G(h)$ has an arc from $i$ to every $j\in [n]$. Thus, $G(h)=K_n$ and
      since $h\sim f$ we have $K_n\in\G(f)$.

    \item Case 2: $f(a)=a$ and $f(b)=b$. Let $a^1,\dots,a^n\in A$, all
      distinct. Then $a^1,\dots,a^n,a,b$ are all distinct since $b=f(b)\neq a$.
      So there is a permutation $\pi$ of $\Q{q}^n$ such that:

      \begin{itemize}
        \item $\pi(a)=\ONE$,
        \item $\pi(b)=\ZERO$,
        \item $\pi(a^i)=e_i$ for $1\leq i\leq n$.
      \end{itemize}
      See Figure \ref{fig:independent_set_1_case2} for an illustration. Let
      $h=\pi\circ f\circ\pi^{-1}$ and let us prove that $G(h)=K_n$. For $i\in
      [n]$, $h(\ZERO)=\ZERO$ and $h(e_i)=\ONE$ differ in every component, hence
      $G(h)$ has an arc from $i$ to every $j\in [n]$. Thus, $G(h)=K_n$ and
      since $h\sim f$ we have $K_n\in\G(f)$.

    \item Case 3: $f(a)=a$ and $f(b)\neq b$. Since $|A|>n$, there is
      $A'\subseteq A\setminus\{f(b)\}$ of size $n$.  Then $A'\cup\{b\}$ is an
      independent set of size $n+1$ and $|f(A'\cup\{b\})|=|\{a,f(b)\}|=2$, thus
      $K_n\in \G(f)$ by Lemma~\ref{lem:independent_set_2k}.
  \end{itemize}
\end{proof}

\begin{proof}[\BF{Proof of Lemma~\ref{lem:independent_set}}]
  Suppose that $f$ has an independent set $A$ of size $|A|\geq 2n$. Then we can
  choose $A$ so that $|A|=2n$, and then $|f(A)|\leq 2n$. Consequently, if
  $|f(A)|$ is even, we have $|f(A)|=2k$ for some $k\in [n]$ and thus
  $K_n\in\G(f)$ by Lemma~\ref{lem:independent_set_2k}. So assume that $|f(A)|$
  is odd. If $|f(A)|=1$ then $K_n\in\G(f)$ by
  Lemma~\ref{lem:independent_set_1}. So assume that $|f(A)|=2k+1$ for some
  $1\leq k<n$. Let us write $f(A)=\{a^1,\dots,a^{2k+1}\}$, and let
  $A_p=f^{-1}(a^p)\cap A$ for $1\leq p\leq 2k+1$. Suppose, without loss, that
  $|A_1|\leq |A_2|\leq \dots\leq |A_{2k+1}|$. Then $A'=A\setminus A_1$ is an
  independent set with $|f(A')|=2k$. If $|A'|\geq n+k$ then $K_n\in\G(f)$ by
  Lemma~\ref{lem:independent_set_2k}. So suppose that $|A'|<n+k$. Since
  $|A'|=2n-|A_1|$ we have $|A_1|> n-k$. If $|A_2|\geq k$ then $A''=A_1\cup A_2$
  is an independent with $|A''|\geq n+1$ and $|f(A'')|=2$ thus $K_n\in\G(f)$ by
  Lemma~\ref{lem:independent_set_2k}. So assume that $|A_2|<k$. Then $n-k<
  |A_1|\leq|A_2|< k$, so $2k>n$. We deduce that  
  \[
    2n=|A|=|A_1|+\cdots+|A_{2k+1}|\geq (2k+1)|A_1|> n|A_1|.
  \]
  So $|A_1|=1$ but then $|A'|=2n-1\geq n+k$, a contradiction. Thus,
  $K_n\in\G(f)$ in every case. 
\end{proof}

By Lemmas~\ref{lem:fixed_points}, \ref{lem:limit_cycles} and
\ref{lem:independent_set}, to prove Theorem~\ref{thm:complete} for $n\geq 5$ it
is sufficient to prove the following. 

\begin{lemma}
  Let $f\in F(n,q)$ with $f\neq\cst,\id$ and $n\geq 5$. Then at least one of
  the following holds:
  \begin{itemize}
    \item $f$ has at least $2n$ fixed points,
    \item $f$ has at least $n$ limit cycles of length $\geq 3$, 
    \item $f$ has an independent set of size at least $2n$.
  \end{itemize}
\end{lemma}

\begin{proof}
  Let $F$ be the set of fixed points of $f$, and let $L$ be a minimal subset
  of $\Q{q}^n$ intersecting every limit cycle of $f$ of length $\geq 3$. Let
  $\Gamma'$ be obtained from $\Gamma(f)$ by deleting the vertices in $F\cup L$;
  then $\Gamma'$ has only cycles of length two, thus it is bipartite. Suppose
  that the first two conditions are false, that is, $|F|<2n$ and $|L|<n$. Then
  $\Gamma'$ has at least $N=q^n-3n+2$ vertices. Since $\Gamma'$ is bipartite,
  it has an independent set $A$ of size at least $N/2$. Since $\Gamma'$ is an
  induced subgraph of $\Gamma(f)$, $A$ is an independent set of $f$, and since
  $n\geq 5$ we have 
  \[
    |A|\geq \lceil N/2 \rceil\geq \lceil(2^n-3n+2)/2\rceil\geq 2n
  \]
  thus the third condition holds.   
\end{proof}

\section{Minimum in-degree}\label{sec:in-degree}

In this section, we will prove Theorem \ref{thm:in-degree}, that we restate.
Given $f\in F(n,q)$, we denote by \EM{$\delta^-(f)$} the minimum integer $d$
such that some digraph in $\G(f)$ has minimum in-degree at most $d$. The
quantity \EM{$\delta^-(n,q)$} defined in the introduction is then the maximum
of $\delta^-(f)$ for $f\in F(n,q)$, and Theorem~\ref{thm:in-degree} says that,
for any $n,q\geq 2$ with $n\geq 3$ or $q\geq 3$, we have 
\[
  2\leq\delta^-(n,2)\leq 5
  \qquad\textrm{and}\qquad
  2\leq\delta^-(n,q)\leq (2+o(1))q.
\] 

\medskip
Instead of considering the in-degree, it is more convenient to consider the
strict in-degree, defined as follows. Given a digraph $G$ on $[n]$, the
\EM{strict in-degree} of $i\in [n]$ in $G$ is the number of in-neighbors of $i$
in $G$ {\em distinct} from $i$. The minimum strict in-degree of $G$ is then the
minimum strict in-degree of a vertex in $G$. Given $f\in F(n,q)$, we denote by
\EM{$\delta^-_s(f)$} the smallest integer $d$ such that some digraph in $\G(f)$
has minimum strict in-degree at most $d$; we have $\delta^-(f)-1\leq
\delta^-_s(f)\leq \delta^-(f)$. 

\medskip
Our first goal is to give, for any $0\leq k<n$, a necessary and sufficient
condition for $\delta^-_s(f)\leq n-k-1$ which is invariant by isomorphism. In
this way, we have a description of $\delta^-_s(f)$ which only relies on the
isomorphic class of $f$. This description is based on the following definition. 

\begin{definition}[$k$-nice partition]
  Let $f \in F(n,q)$ and $k\in [n]$. A \EM{$k$-nice partition} of $f$ is a
  partition $A=\{A_1,\dots,A_q\}$ of $\Q{q}^n$, with some parts possibly empty,
  such that $|A_p|$ and $|A_\ell\cap f^{-1}(A_p)|$ are multiples of $q^k$ for
  all $p,\ell\in[q]$. 
\end{definition}

A partition is \EM{balanced} when its parts have all the same size. The
characterization is the following. 

\begin{lemma} \label{lem:nice}
  For every $f \in F(n,q)$ and $0\leq k< n$, we have $\delta^-_s(f)\leq n-k-1$
  if and only if $f$ has a balanced $k$-nice partition.
\end{lemma}

Actually, the proof of Theorem~\ref{thm:in-degree} only uses one direction,
that $\delta^-_s(f)\leq n-k-1$ if $f$ has a balanced $k$-nice partition. But
the main property needed for the other direction will be used later, and so we
give it separately (in Lemma~\ref{lem:nice_1} below), and then put things
together to get the characterization. For that we use the following
definitions. Given $x,y\in\Q{q}^n$, we denote by \EM{$\Delta(x,y)$} the set of
$i\in [n]$ with $x_i\neq y_i$. Given $I\subseteq [n]$, a set $X \subseteq
\Q{q}^n$ is \EM{$I$-closed} if, for any $x,y\in \Q{q}^n$ with
$\Delta(x,y)\subseteq I$, we have $x\in X$ if and only if $y\in X$. It is easy
to see that if $X$ is $I$-closed then $|X|$ is a multiple of $q^{|I|}$, and
that if $X$ and $Y$ are both $I$-closed then so is $X \cap Y$.

\begin{lemma}\label{lem:nice_1}
  Let $f \in F(n,q)$ and suppose that $G(f)$ has a vertex $i$ with strict
  in-degree at most $n-k-1$. For $1\leq p\leq q$, let $A_p$ be the set of
  $x\in\Q{q}^n$ with $x_i = p-1$. Then $A=\{A_1,\dots,A_q\}$ is a balanced
  $k$-nice partition of $f$.
\end{lemma}

\begin{proof}
  For $k=0$ the result is obvious, so assume that $k\geq 1$. Since $A$ is
  obviously balanced, we just have to prove that $A$ is $k$-nice. Since $i$ is
  of strict in-degree at most $n-k-1$, there is a set $I\subseteq [n]$ of size
  $k$ such that $G(f)$ has no arc from $I$ to $i$ and $i\notin I$. Since
  $i\notin I$, each $A_p$ is $I$-closed. Let $x$ be any configuration in
  $f^{-1}(A_p)$, that is, $f_i(x)=p-1$. Since there are no arc from $I$ to $i$,
  for every $y\in\Q{q}^n$ with $\Delta(x,y)\subseteq I$ we have
  $f_i(y)=f_i(x)=p-1$ and thus $y\in f^{-1}(A_p)$. As such, $f^{-1}(A_p)$ is
  $I$-closed. Since $A_\ell$ and $f^{-1}(A_p)$ are  $I$-closed for any
  $p,\ell\in [k]$, their intersection $A_\ell \cap f^{-1}(A_p)$ is also
  $I$-closed, and so $|A_\ell\cap f^{-1}(A_p)|$ is a multiple of $q^k$, as
  desired. 
\end{proof}

\begin{proof}[\BF{Proof of Lemma~\ref{lem:nice}}]
  For $k=0$ the result is obvious, so assume that $k\geq 1$.

  \medskip
  First, let $f \in F(n,q)$ and suppose that the minimum strict in-degree of
  $G(h)$ is at most $n-k-1$ for some $h\sim f$. By Lemma~\ref{lem:nice_1}, $h$
  has a balanced $k$-nice partition and this trivially implies that $f$ has a
  balanced $k$-nice partition. This proves the first direction. 


  \medskip
  For the other direction, let $f \in F(n,q)$ and suppose that $f$ has a
  balanced $k$-nice partition $A = \{A_1,\dots,A_q\}$. For all $p,\ell\in [q]$,
  let $A^-_p = f^{-1}(A_p)$ and $a_{p,\ell} = |A_p \cap A^-_\ell| / q^k$, which
  is an integer since $|A_p \cap A^-_\ell|$ is a multiple of $q^k$. Note that,
  for all $p\in [q]$, we have 
  \[
    \sum\limits_{\ell\in [q]}a_{p,\ell} = q^{-k}\sum_{\ell\in [q]}|A_p \cap A^-_\ell|=q^{-k}|A_p| = q^{n-k-1}.
  \]
  Hence, there is a partition $X_p = \{X_{p,1}, \dots, X_{p,q}\}$ of
  $\Q{q}^{n-k-1}$ (with some members possibly empty) such that $|X_{p,\ell}| =
  a_{p,\ell}$ for all $\ell\in [q]$. For all $p,\ell\in [q]$, let 
  \[
    Y_{p,\ell} = \{x\in\Q{q}^n \mid x_n=p-1,~(x_1,\dots,x_{n-k-1}) \in X_{p,\ell}\}.
  \]
  Since each configuration in $X_{p,\ell}$ is extended into exactly $q^k$
  configurations in $Y_{p,\ell}$, we have $|Y_{p,\ell}|=a_{p,\ell}
  q^k=|A_p \cap A_\ell^-|$. Furthermore, it is clear that $Y_{p,\ell}\cap
  Y_{p',\ell'}=\emptyset$ for all $p,\ell,p',\ell'\in [q]$ with $(p,\ell) \neq
  (p',\ell')$, so $\{Y_{p,\ell}\mid p,\ell\in [q]\}$ is a partition of $\Q{q}^n$.
  Hence, there is a permutation $\pi$ of $\Q{q}^n$ such that $\pi(A_p \cap
  A_\ell^-) = Y_{p,\ell}$ for all $p,\ell\in [q]$. 

  \medskip
  Let $h=\pi\circ f\circ\pi^{-1}$ and let us prove that $G(h)$ has no arc from
  $I=[n-k,n-1]$ to $n$. So consider any $x,y\in\Q{q}^n$ with
  $\Delta(x,y)\subseteq I$, and let us prove that $h_n(x)=h_n(y)$. Let
  $p,\ell\in [q]$ such that $x \in Y_{p,\ell}$. Since $Y_{p,\ell}$ is
  $I$-closed, we also have $y\in Y_{p,\ell}$. Hence, setting $x'=\pi^{-1}(x)$
  and $y'=\pi^{-1}(y)$, we have $x',y'\in A_p\cap A^-_\ell$, and thus
  $f(x'),f(y')\in A_\ell$. So there are $r,s\in [q]$ such that $f(x')\in
  A_\ell\cap A^-_r$ and $f(y')\in A_\ell\cap A^-_s$, and we deduce that
  $h(x)=\pi(f(x'))\in Y_{\ell,r}$ and $h(y)=\pi(f(x'))\in Y_{\ell,s}$.
  Consequently, $h_n(x)=h_n(y)=\ell-1$ as desired. Hence, $G(h)$ has no arc
  from $I$ to $n$, and thus the strict in-degree of $n$ in $G(h)$ is at most
  $n-k-1$. 
\end{proof}

\subsection{Upper bounds}

By the preceding characterization, to prove that $\delta^-_s(f)$ is small we
have to prove that $f$ has a balanced $k$-nice partition for some large $k$. We
will prove that using two tools. The first is the following easy lemma. Given
$f \in F(n,q)$, we say that $X \subseteq \Q{q}^n$ is a \EM{weakly independent
set} of $f$ if $f(x)\notin X\setminus\{x\}$ for all $x\in X$ (so $X$ is an
independent set that can contain loops). 

\begin{lemma}\label{prop:independant}
  Let $f \in F(n,q)$ and $Y \subseteq \Q{q}^n$. Then $f$ has a weakly
  independent set $X\subseteq Y$ with $|X| \geq |Y|/3$.
\end{lemma}

\begin{proof}
  For $X\subseteq \Q{q}^n$, we denote by $\Gamma[X]$ the subgraph of
  $\Gamma(f)$ induced by $X$. Let $Z$ be a minimal subset of $Y$ intersecting
  every cycle of $\Gamma[Y]$ of length at least $3$. Let $Y' = Y \setminus Z$.
  Since $|Z| \leq |Y|/3$ we have $|Y'| = |Y| - |Z| \geq 2|Y|/3$. Let $\Gamma'$
  be obtained from $\Gamma[Y']$ by removing loops. Then $\Gamma'$  is
  bipartite, as every cycle is of length~$2$. So $\Gamma'$ has an independent
  set $X \subseteq Y'$ with $|X| \geq |Y'|/2 \geq |Y|/3$. Since $X$ is an
  independent set of $\Gamma'$, there are no distinct $x,y\in X$ such that
  $f(x)=y$ and we deduce that $X$ is a weakly independent set of $f$.
\end{proof}

The second tool comes from additive number theory and is based on the following definition. 

\begin{definition}[$s(m,d)$]
  For any integers $m,d\geq 1$, let \EM{$s(m,d)$} be the smallest integer $s$
  such that, for any functions $a_1,\dots,a_d:[s]\to\mathbb{N}$, there is a
  subset $X\subseteq [s]$ of size $m$ such that $\sum_{x\in X}a_i(x)$ is a
  multiple of $m$ for every $i\in [d]$.
\end{definition}

That $s(m,d)$ exists, and is at most $(m-1)m^d-1$, is an easy consequence of
the pigeonhole principle. Indeed, let us say that $x,y\in [s]$ are equivalent
if, for all $i\in [d]$, we have $a_i(x)\equiv a_i(y)\mod m$. There are at most
$m^d$ equivalence classes and so if $s\geq (m-1)m^d+1$ then, by the pigeonhole
principle, some class is of size at least $m$, and any subset $X$ of size $m$
in this class is as in the definition. But much better bounds exist. Erd\"{o}s,
Ginzburg and Ziv famously proved that $s(m,1)=2m-1$ \cite{EGZ61}. The exact
value of $s(m,2)$ has been established by Reiher, resolving a longstanding
conjecture of Kemnitz.

\begin{theorem}[Reiher \cite{Reiher07}] \label{thm:Reiher}
  For every $m\geq 1$, 
  \[
    s(m,2)=4m-3. 
  \]
\end{theorem}

Hence, $s(m,1)$ and $s(m,2)$ are linear in $m$. Using deep arguments, Alon and
Dubiner proved that $s(m,d)$ is actually linear in $m$ for every fixed $d$. 

\begin{theorem}[Alon-Dubiner \cite{AD95}] \label{thm:Alon-Dubiner}
  There exists an absolute constant $c$ such that, for all $m,d\geq 1$, 
  \[
    s(m,d)\leq (cd\log_2d)^dm.
  \]
\end{theorem}

We now prove that any $f\in F(n,q)$ has a balanced $k$-nice partition for some
large $k$, expressed with the quantity $s(m,d)$. We begin with the case $q=2$,
where we optimize the use of $s(m,d)$ in order to use the exact value given by
Reiher for $d=2$. The general case, where such optimization is not needed,
follows the same approach and is presented afterwards.

\begin{lemma}\label{lem:alwaysnice2}
  Every $f \in F(n,2)$ has a balanced $k$-nice partition for any $k$ such that 
  \[
    3s(2^k,2)\leq 2^{n-1}.
  \]
\end{lemma}

\begin{proof}
  Let $f\in F(n,2)$. For $X\subseteq \B^n$, we use $X^-$ as a shorthand for
  $f^{-1}(X)$. Let $k$ be as in the statement. Observe that, by setting
  $A_1=\B^n$ and $A_2=\emptyset$, we obtain a $k$-nice partition $\{A_1,A_2\}$
  since $|A_1|=|A_1\cap A_1^-|=2^n$ and $|A_2|=|A_1\cap A_2^-|=|A_2\cap
  A_1^-|=|A_2\cap A_2^-|=0$. So there is a $k$-nice partition $\{A_1,A_2\}$ with
  $|A_1|\geq 2^{n-1}$ where $A_1$ is minimal for this property. We will prove
  that $|A_1|=2^{n-1}$ so that the partition is balanced. Suppose, for a
  contradiction, that $|A_1|>2^{n-1}$.

  \medskip
  By Lemma \ref{prop:independant}, there is a subset $Y\subseteq A_1$ of size
  at least $|A_1|/3$ which is a weakly independent set. For $x\in Y$, let 
  \begin{align*}
    \alpha_1(x)&=|A_1\cap f^{-1}(x)|,\\[2mm]
    \alpha_2(x)&=|A_2\cap f^{-1}(x)|,\\[2mm]
    \beta_1(x)&=
    \left\{
      \begin{array}{ll}
        1&\textrm{if $f(x)\in A_1\setminus \{x\}$,}\\[1mm]
        0&\textrm{otherwise,}
      \end{array}
      \right.
  \end{align*}

  \medskip
  Since $|A_1|\geq 2^{n-1}$ we have $|Y|\geq |A_1|/3\geq 2^{n-1}/3\geq
  s(2^k,2)$. Considering the $2$ functions $a_1,a_2$ defined on $Y$ by
  $a_1(x)=\alpha_1(x)+\beta_1(x)$ and $a_2(x)=\alpha_2(x)-\beta_1(x)$, we
  deduce from the definition of $s$ that there is a subset $X \subseteq Y$ of
  size $2^k$ such that $\sum_{x\in X} a_1(x)$ and $\sum_{x\in X} a_2(x)$ are
  multiples of $2^k$. Setting $\alpha_i(X)=\sum_{x\in X}\alpha_i(x)$ for
  $i=1,2$ and $\beta_1(X)=\sum_{x\in X}\beta_1(x)$ we obtain that
  $\alpha_1(X)+\beta_1(X)$ and $\alpha_2(X)-\beta_1(X)$ are multiples of $2^k$. 
  Furthermore:
  \begin{align}
    \alpha_i(X)&=|A_i\cap X^-|\textrm{ for $i=1,2$,}\\[1mm]
    \beta_1(X)&=|X\cap A_1^-|-|X\cap X^-|.
  \end{align}
  (1) is obvious, since 
  \[
    |A_i\cap X^-|=|A_i\cap(\cup_{x\in X} f^{-1}(x))|=\sum_{x\in X}|A_i\cap f^{-1}(x)|=\alpha_i(X).
  \]
  For (2), note that $\beta_1(X)$ is the number of $x\in X$ with $f(x)\in
  A_1\setminus \{x\}$. Since $X$ is weakly independent, $f(x)\in A_1\setminus
  \{x\}$ is equivalent to $f(x)\in A_1\setminus X$, and thus $\beta_1(X)$ is
  the number of $x\in X$ with $f(x)\in A_1\setminus X$, that is,
  $\beta_1(X)=|X\cap (A_1\setminus X)^-|=|X\cap A_1^-|-|X\cap X^-|$. 

  \medskip
  Consider the partition
  \[
    \{B_1,B_2\}=\{A_1\setminus X,A_2\cup X\}.
  \]
  Recall that $\{A_1,A_2\}$ is $k$-nice, hence $|A_1|$ is a multiple of $2^k$,
  and since we also have $|A_1|>2^{n-1}$, it follows that
  $|B_1|=|A_1|-2^k\geq 2^{n-1}$. 
  As a consequence, in order to obtain the desired
  contradiction, it is sufficient to prove that $\{B_1,B_2\}$ is $k$-nice. First,
  since $|X|=2^k$, it is clear that $|B_1|$ and $|B_2|$ are multiples of $2^k$.
  So it remains to prove that $|B_i \cap B_j^-|$ is a multiple of $2^k$ for
  every $1\leq i,j\leq 2$. First, using (1) and (2), we have:
  \begin{align*}
    |B_1\cap B_1^-|&=|(A_1\setminus X)\cap (A_1^-\setminus X^-)|\\
    &=|A_1\cap A_1^-|-|A_1\cap X^-|-|X\cap A_1^-|+|X\cap X^-|\\
    &=|A_1\cap A_1^-|-\alpha_1(X)-\beta_1(X)\\[2mm]
    |B_2\cap B_1^-|&=|(A_2\cup X)\cap (A_1^-\setminus X^-)|\\
    &= |A_2\cap A_1^-|-|A_2\cap X^-|+|X\cap A_1^-|-|X\cap X^-|\\
    &= |A_2\cap A_1^-|-\alpha_2(X)+\beta_1(X).
  \end{align*}
  Since $|A_1\cap A_1^-|$, $|A_2\cap A_1^-|$, $\alpha_1(X)+\beta_1(X)$ and
  $\alpha_2(X)-\beta_1(X)$ are multiples of $2^k$, we deduce that $|B_1\cap
  B_1^-|$ and $|B_2\cap B_1^-|$ are multiples of $2^k$. Then we have 
  \begin{align*}
    |B_1\cap B_2^-|&=|B_1|-|B_1\cap B_1^-|\\[2mm]
    |B_2\cap B_2^-|&=|B_2|-|B_2\cap B_1^-|.
  \end{align*}
  Since $|B_1|$, $|B_2|$, $|B_1\cap B_1^-|$ and $|B_2\cap B_1^-|$ are
  multiples of $2^k$, we deduce that $|B_1\cap B_2^-|$ and $|B_2\cap B_2^-|$
  are multiples of $2^k$. Thus, $\{B_1,B_2\}$ is $k$-nice, a contradiction. 
\end{proof}

Combining this lemma with Lemma~\ref{lem:nice} and Reiher's theorem, which says
that $s(m,2)=4m-3$, we obtain the first upper bound in
Theorem~\ref{thm:in-degree}.

\begin{lemma}
  For all $n\geq 1$ we have $\delta^-(n,2)\leq 5$, and $\delta^-(5,2)\leq 4$.
\end{lemma}

\begin{proof}
  Let $f\in F(n,2)$. Using Theorem~\ref{thm:Reiher} we have 
  \[
    3s(2^{n-5},2)=12\cdot2^{n-5}-9\leq 16\cdot 2^{n-5}=2^{n-1}. 
  \]
  Hence, by Lemma~\ref{lem:alwaysnice2}, $f$ has a balanced $(n-5)$-nice
  partition. Thus, $\delta^-_s(f)\leq n-(n-5)-1=4$ by Lemma~\ref{lem:nice}.
  This proves that $\delta^-(n,2)\leq 5$. Then, since $3s(2,2)=15\leq 2^4$, we
  deduce from Lemma~\ref{lem:alwaysnice2} that every $f\in F(5,2)$ has a
  $1$-nice partition and thus $\delta^-_s(f)\leq 3$ by Lemma~\ref{lem:nice}.
  This proves that $\delta^-(5,2)\leq 4$. 
\end{proof}

We now extend the previous arguments to larger alphabets.

\begin{lemma}\label{lem:alwaysniceq3}
  Every $f \in F(n,q)$ has a balanced $k$-nice partition for any $k$ such that 
  \[
    3 s(q^k,2q)\leq q^{n-1}.
  \]
\end{lemma}

\begin{proof}
  Let $f\in F(n,q)$. For $X\subseteq \Q{q}^n$, we use $X^-$ as a shorthand for
  $f^{-1}(X)$. Let $k$ be as in the statement. Given a partition
  $A=\{A_1,\dots,A_q\}$ of $\Q{q}^n$, we say that $A$ is {\em $1$-dominating} if
  $|A_1|\geq q^{n-1}$ and $|A_i|\leq q^{n-1}$ for all $2\leq i\leq q$. Observe
  that, by setting $A_1=\Q{q}^n$ and $A_i=\emptyset$ for $2\leq i\leq q$, we
  obtain a $1$-dominating $k$-nice partition: for every $2\leq i\leq q$ and
  $1\leq j\leq q$, we have $|A_1|=|A_1\cap A_1^-|=q^n$ and $|A_i|=|A_1\cap
  A_i^-|=|A_i\cap A_j^-|=0$. So we can suppose that $A$ is a $1$-dominating
  $k$-nice partition with $|A_1|$ minimal. 

  \medskip
  We will prove that $|A_1|=q^{n-1}$, which (since $A$ is $1$-dominating)
  forces $A$ to be balanced. Suppose, for a contradiction, that
  $|A_1|>q^{n-1}$.  Since $A$ is $1$-dominating there is $2\leq i\leq n$ such
  that $|A_i|<q^{n-1}$. Suppose, without loss, that $|A_2|<q^{n-1}$. 

  \medskip
  By Lemma~\ref{prop:independant}, there is a subset $Y\subseteq A_1$ of
  size at least $|A_1|/3$ which is a weakly independent set. For $x\in Y$ and
  $1\leq i\leq q$ let 
  \begin{align*}
    \alpha_i(x)&=|A_i\cap f^{-1}(x)|,\\[2mm]
    \beta_i(x)&=
    \left\{
      \begin{array}{ll}
        1&\textrm{if $f(x)\in A_i\setminus \{x\}$,}\\[1mm]
        0&\textrm{otherwise,}
      \end{array}
      \right.\\[2mm]
      \beta'_i(x)&=
      \left\{
        \begin{array}{ll}
          \beta_i(x)+1&\textrm{if $i=2$ and $f(x)=x$,}\\[1mm]
          \beta_i(x)&\textrm{otherwise.}
        \end{array}
        \right.
  \end{align*}

  \medskip
  Since $|A_1|\geq q^{n-1}$ we have $|Y|\geq |A_1|/3\geq q^{n-1}/3\geq
  s(q^k,2q)$. Considering the $2q$ functions $\alpha_i,\beta'_i$, we deduce
  from the definition of $s$ that there is a subset $X \subseteq Y$ of size
  $q^k$ such that $\alpha_i(X)=\sum_{x\in X}\alpha_i(x)$ and
  $\beta'_i(X)=\sum_{x\in X}\beta'_i(x)$ are multiples of $q^k$ for all $1\leq
  i\leq q$. Furthermore:
  \begin{align}
    \alpha_i(X)&=|A_i\cap X^-|\textrm{ for all $1\leq i\leq q$,}\tag{1}\\[1mm]
    \beta_i(X)&=|X\cap A_i^-|\textrm{ for all $2\leq i\leq q$,}\tag{2}\\[1mm] 
    \beta'_i(X)&=|X\cap A_i^-|\textrm{ for all $3\leq i\leq q$,}\tag{3}\\[1mm]
    \beta'_1(X)&=|X\cap A_1^-|-|X\cap X^-|,\tag{4}\\[1mm]
    \beta'_2(X)&=|X\cap A_2^-|+|X\cap X^-|.\tag{5}
  \end{align}
  (1) is obvious, since 
  \[
    |A_i\cap X^-|=|A_i\cap(\cup_{x\in X} f^{-1}(x))|=\sum_{x\in X}|A_i\cap f^{-1}(x)|=\alpha_i(X).
  \]
  For (2), note that $\beta_i(X)$ is the number of $x\in X$ with $f(x)\in
  A_i\setminus\{x\}$, which is equivalent to $f(x)\in A_i$ since $2\leq i\leq
  q$ and $x\in A_1$. So $\beta_i(X)$ is the number of $x\in X$ with $f(x)\in
  A_i$, that is, $\beta_i(X)=|X\cap A_i^-|$. (3) results from (2) since
  $\beta'_i=\beta_i$ for $3\leq i\leq q$. For (4), note that $\beta'_1=\beta_1$
  and that $\beta_1(X)$ is the number of $x\in X$ with $f(x)\in A_1\setminus
  \{x\}$. Since $X$  is weakly independent, $f(x)\in A_1\setminus \{x\}$ is
  equivalent to $f(x)\in A_1\setminus X$, and thus $\beta'_1(X)$ is the number
  of $x\in X$ with $f(x)\in A_1\setminus X$, that is, $\beta'_1(X)=|X\cap
  (A_1\setminus X)^-|=|X\cap A_1^-|-|X\cap X^-|$. To prove (5), let $\ell$ be the
  number of $x\in X$ with $f(x)=x$. It is clear that
  $\beta'_2(X)=\beta_2(X)+\ell$, and using the second point we obtain
  $\beta'_2(X)=|X\cap A_2^-|+\ell$. But since $X$ is weakly independent, for $x\in
  X$ we have $f(x)=x$ if and only if $f(x)\in X$, and thus $\ell=|X\cap X^-|$ so
  that $\beta'_2(X)=|X\cap A_2^-|+|X\cap X^-|$. 

  \medskip
  Consider the partition
  \[
    B=\{B_1,B_2,B_3,\dots,B_q\}=\{A_1\setminus X,A_2\cup X,A_3,\dots,A_q\}.
  \]
  Recall that $A$ is $k$-nice, hence $|A_1|,|A_2|$ are multiples of $q^k$,
  and since we also have $|A_1|>q^{n-1}$ and $|A_2|<q^{n-1}$,
  it follows that $|B_1|=|A_1|-q^k\geq q^{n-1}$ and $|B_2|=|A_2|+q^k\leq q^{n-1}$,
  therefore $B$ is $1$-dominating. As a consequence, in order to
  obtain the desired contradiction, it is sufficient to prove that $B$ is
  $k$-nice. First, since $|X|=q^k$, it is clear that $|B_i|$ is a multiple of
  $q^k$ for every $1\leq i\leq q$. So it remains to prove that $|B_i \cap
  B_j^-|$ is a multiple of $q^k$ for every $1\leq i,j\leq q$. This is obvious
  for $3\leq i,j\leq q$ since $|B_i \cap B_j^-|=|A_i \cap A_j^-|$. We show that
  the other intersection sizes are multiples of $q^k$ by expressing these sizes
  as a sum of (positive or negative) numbers that are already known to be
  multiples of $q^k$. For that we use many times the properties
  (1),(3),(4),(5). We begin with the case $1\leq i,j\leq 2$:
  \begin{align*}
    |B_1\cap B_1^-|&=|(A_1\setminus X)\cap (A_1^-\setminus X^-)|\\
    &=|A_1\cap A_1^-|-|A_1\cap X^-|-|X\cap A_1^-|+|X\cap X^-|\\
    &=|A_1\cap A_1^-|-\alpha_1(X)-\beta'_1(X)\\[2mm]
    |B_1\cap B_2^-|&=|(A_1\setminus X)\cap (A_2^-\cup X^-)|\\
    &=|A_1\cap A_2^-|+|A_1\cap X^-|-|X\cap A_2^-|-|X\cap X^-|\\
    &=|A_1\cap A_2^-|+\alpha_1(X)-\beta'_2(X)\\[2mm]
    |B_2\cap B_1^-|&=|(A_2\cup X)\cap (A_1^-\setminus X^-)|\\
    &= |A_2\cap A_1^-|-|A_2\cap X^-|+|X\cap A_1^-|-|X\cap X^-|\\
    &= |A_2\cap A_1^-|-\alpha_2(X)+\beta'_1(X)\\[2mm]
    |B_2\cap B_2^-|&=|(A_2\cup X)\cap (A_2^-\cup X^-)|\\
    &= |A_2\cap A_2^-|+|A_2\cap X^-|+|X\cap A_2^-|+|X\cap X^-|\\
    &= |A_2\cap A_2^-|+\alpha_2(X)+\beta'_2(X).
  \end{align*}
  For $1\leq i\leq 2$ and $3\leq j\leq q$ we have:
  \begin{align*}
    |B_1\cap B_j^-|&=|(A_1\setminus X)\cap A_j^-|\\
    &= |A_1\cap A_j^-|-|X\cap A_j^-|\\
    &= |A_1\cap A_j^-|-\beta'_j(X)\\[2mm]
    |B_2\cap B_j^-|&=|(A_2\cup X)\cap A_j^-|\\
    &= |A_2\cap A_j^-|+ |X\cap A_j^-|\\
    &= |A_2\cap A_j^-|+\beta'_j(X).
  \end{align*}
  Finally, for $3\leq i\leq q$ and $1\leq j\leq 2$, we have:
  \begin{align*}
    |B_i\cap B_1^-|&=|A_i\cap (A_1^-\setminus X^-)|\\
    &= |A_i\cap A_1^-|- |A_i\cap X^-|\\
    &= |A_i\cap A_1^-|- \alpha_i(X)\\[2mm]
    |B_i\cap B_2^-|&=|A_i\cap (A_2^-\cup X^-)|\\
    &= |A_i\cap A_2^-|+ |A_i\cap X^-|\\
    &= |A_i\cap A_2^-|+ \alpha_i(X).
  \end{align*}
  So $|B_i\cap B_j^-|$ is indeed a multiple of $q^k$ for all $1\leq i,j\leq q$,
  thus $B$ is $k$-nice, a contradiction. 
\end{proof}

Combining this lemma with Lemma~\ref{lem:nice} and Alon-Dubiner's theorem,
which says that $s(m,d)\leq (cd\log_2d)^dm$ for some constant $c$, we obtain
the second upper bound in Theorem~\ref{thm:in-degree}.

\begin{lemma}
  For all $n\geq 1$ and $q\geq 2$, we have $\delta^-(n,q)\leq (2+o(1))q$.
\end{lemma}

\begin{proof}
  Let $f\in F(n,q)$, let $\alpha(q)=2c\log_2(2q)$ where $c$ is as in Theorem
  \ref{thm:Alon-Dubiner}, and let $k=\lfloor \ell\rfloor$ with 
  \[
    \ell=n-2q\log_q(q\alpha(q))-1-\log_q3.
  \] 
  Using Theorem \ref{thm:Alon-Dubiner} for the first inequality, we have 
  \[
    3s(q^k,2q)\leq 3(q\alpha(q))^{2q} q^k\leq 3(q\alpha(q))^{2q} q^\ell=q^{n-1}.
  \]
  Hence, by Lemma~\ref{lem:alwaysniceq3}, $f$ has a balanced $k$-nice
  partition, and thus $\delta^-_s(f)\leq n-k-1$ by Lemma~\ref{lem:nice}. So
  $\delta^-(f)\leq n-k$, and since $k\geq \ell-1$ we have
  \[
    \delta^-(f)\leq  2q\log_q(q\alpha(q))+2+\log_q3=\left(2+2\log_q(\alpha(q))+(2+\log_q3)q^{-1}\right) q.
  \]
  Since $2\log_q(\alpha(q))+(2+\log_q3)q^{-1}$ tends to $0$ as $q$ tends to
  infinity, this proves the lemma. 
\end{proof}

\subsection{Lower bound}

It remains to prove the lower bound in Theorem~\ref{thm:in-degree}, which is
based on a simple construction. 

\begin{lemma}\label{lem:lower_bound}
  For all $n,q\geq 2$ with $n\geq 3$ or $q\geq 3$, we have $\delta^-(n,q)\geq 2$. 
\end{lemma}

\begin{figure}[h]
  \centering
  \begin{tikzpicture}
    \node[circle,inner sep=1] (x1) at ({-90-0*72}:1.2){$x^1$};
    \node[circle,inner sep=1] (x2) at ({-90-1*72}:1.2){$x^2$};
    \node[circle,inner sep=1] (x3) at ({-90-2*72}:1.2){$x^3$};
    \node[circle,inner sep=1] (x4) at ({-90-3*72}:1.2){$x^4$};
    \node[circle,inner sep=1] (x5) at ({-90-4*72}:1.2){$x^5$};
    \node (y1) at (0,-2.5){\tiny$\bullet$};
    \node (y2) at (1,-2.5){\tiny$\bullet$};
    \node (y3) at (-1,-2.5){\tiny$\bullet$};
    \path[thick,->,draw,black]
    (x1) edge (x2)
    (x2) edge (x3)
    (x3) edge (x4)
    (x4) edge (x5)
    (x5) edge (x1)
    (y1) edge (x1)
    (y2) edge (x1)
    (y3) edge (x1)
    ;
  \end{tikzpicture}
  \caption{\label{fig:degmin_2more}
  The function of Lemma \ref{lem:lower_bound} for $n = 3$ and $q = 2$.
  }
\end{figure}
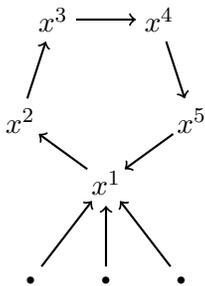

\begin{proof}
  Let $f\in F(n,q)$ with $n,q\geq 2$ and $n\geq 3$ or $q\geq 3$. Suppose that
  $f$ has a limit cycle of length $\ell=q^{n-1}+1$, whose configurations are
  $x^1,\dots,x^\ell$ in order, and suppose that $x^1$ is the image of the all
  configurations which are not in the limit cycle; so $f$ has exactly $\ell$
  images. See Figure \ref{fig:degmin_2more} for an illustration. We will prove
  that $G(f)$ has minimum in-degree at least two, and since $f$ is described up
  to isomorphism this proves that $\delta^-(f)\geq 2$ and the lemma follows.
  Suppose, for a contradiction, that $G(f)$ has a vertex $i$ with in-degree at
  most $1$. 

  \medskip
  If $i$ is of in-degree $0$, then $f_i$ is a constant function, which always returns some $a\in\Q{q}$. Hence, $x_i=a$ for all the images $x$ of $f$, and thus $f$ has at most $q^{n-1}<\ell$ images, a contradiction. 

  \medskip
  So $i$ is of in-degree $1$. Let $j$ be its in-neighbor (we may have $i=j$).
  There is then a non-constant function $g:\Q{q}\to\Q{q}$ such that
  $f_i(x)=g(x_j)$ for all $x\in \Q{q}^n$. Let $Y=\{x^1,\dots,x^{\ell-1}\}$ and
  note that $Y=\Q{q}^n\setminus f^{-1}(x^1)$ by construction. For $a\in\Q{q}$
  we set $X_{i,a}=\{x\in\Q{q}^n\mid x_i=a\}$ and $X_{j,a}=\{x\in\Q{q}^n\mid
  x_j=a\}$. If $X_{j,a}$ is not included in $Y$, that is if there is a
  configuration $x$ with $x_j=a$ and $f(x)=x^1$, then $g(a)=f_i(x)=x^1_i$. We
  deduce that if $X_{j,a}$ is not included in $Y$ for every $a\in\Q{q}$, then
  $g$ is a constant function (which always return $x^1_i$), a contradiction. So
  $X_{j,a}\subseteq Y$ for some $a$, and since $Y$ is of size
  $\ell-1=q^{n-1}=|X_{j,a}|$  we deduce that $Y=X_{j,a}$. Hence, setting
  $b=g(a)$, we have $f(Y)\subseteq X_{i,b}$, and since
  $f(Y)=\{x^2,\dots,x^\ell\}$ is of size $\ell-1=q^{n-1}=|X_{i,b}|$ we deduce
  that $f(Y)=X_{i,b}$. Consequently, $X_{j,a}$ and $X_{i,b}$ are distinct.
  Furthermore, they are not disjoint since $X_{j,a}\cap X_{i,b}=Y\cap
  f(Y)=\{x^2,\dots,x^{\ell-1}\}$ is of size $\ell-2=q^{n-1}-1>0$. We deduce
  that $i\neq j$ and thus $|X_{j,a}\cap X_{i,b}|=q^{n-2}$. We obtain
  $q^{n-2}=q^{n-1}-1$, which is a contradiction since $n\geq 3$ or $q\geq 3$.   
\end{proof}

\section{Dense digraphs}\label{sec:dense}

In this section, we will prove Theorem~\ref{thm:dense}, that for $n,q\geq 2$
there is $f\in F(n,q)$ such that every digraph in $\G(f)$ has at least $\lfloor
n^2/4 \rfloor$ arcs. We need two lemmas. In a digraph, the \EM{sources} are the
vertices of in-degree $0$, and the \EM{non-sources} are the vertices of
in-degree at least~$1$. 

\begin{lemma} \label{lem:sum_images}
  Let $f \in F(n,q)$ and let $X$ be the images of $f$. Let $k\in [n]$ and
  suppose that, for all non-empty strict subsets $Y$ of $X$, $|f^{-1}(Y)|$ is
  not a multiple of $q^k$. Then every non-source of $G(f)$ has strict in-degree
  at least $n-k$.
\end{lemma}

\begin{proof}
  Suppose, for a contradiction, that $G(f)$ has a non-source $i$ with strict
  in-degree at most $n-k-1$. For $1\leq p\leq q$, let $A_p$ be the set of
  $x\in\Q{q}^n$ with $x_i = p-1$. By Lemma~\ref{lem:nice_1},
  $A=\{A_1,\dots,A_q\}$ is a balanced $k$-nice partition of $f$. Let $1\leq p\leq
  q$ such that $X\cap A_p\neq \emptyset$ and let $Y=X\cap A_p$. If $X\subseteq
  A_p$, then it means that $f_i$ is a constant function (which always returns
  $p$), which is a contradiction because $i$ has in-degree at least $1$. We
  deduce that $Y$ is a non-empty strict subset of $X$. Furthermore, we have
  $|f^{-1}(Y)|=|f^{-1}(A_p)|=\sum_{\ell=1}^q|A_\ell\cap f^{-1}(A_p)|$. Since
  $A$ is $k$-nice, each term of the sum is a multiple of $q^k$. So
  $|f^{-1}(Y)|$ is a multiple of $q^k$, a contradiction. 
\end{proof}

\begin{lemma}\label{lem:k_non_sources}
  For all $n,q\geq 2$ and $k\in [n]$, there is $f \in F(n,q)$ such that, for
  every $G\in\G(f)$, $G$ has at least $k$ non-sources and all the non-sources
  of $G$ have strict in-degree at least $n-k$.
\end{lemma}

\begin{figure}[h]
  \centering
  \begin{tikzpicture}
    \def\c{2}
    \node[outer sep=1,inner sep=1.5,circle,thick] (x1) at (0,0){$x^*$};
    \node (x2) at (3,0){\tiny$\bullet$};
    \node (x3) at (4,0){\tiny$\bullet$};
    \node (x4) at (5,0){\tiny$\bullet$};
    \node (y1) at (0:\c){\tiny$\bullet$};
    \node (y2) at (20:\c){\tiny$\bullet$};
    \node (y3) at (40:\c){\tiny$\bullet$};
    \node (y4) at (60:\c){\tiny$\bullet$};
    \node (y5) at (80:\c){\tiny$\bullet$};
    \node (y6) at (100:\c){\tiny$\bullet$};
    \node (y7) at (120:\c){\tiny$\bullet$};
    \node (y8) at (140:\c){\tiny$\bullet$};
    \node (y9) at (160:\c){\tiny$\bullet$};
    \node (y10) at (180:\c){\tiny$\bullet$};
    \node (y11) at (200:\c){\tiny$\bullet$};
    \node (y12) at (-20:\c){\tiny$\bullet$};
    \draw[->,thick] (x1.-112) .. controls (-0.5,-0.8) and (0.5,-0.8) .. (x1.-68);
    \draw[->,thick] (x2.-112) .. controls (2.5,-0.7) and (3.5,-0.7) .. (x2.-68);
    \draw[->,thick] (x3.-112) .. controls (3.5,-0.7) and (4.5,-0.7) .. (x3.-68);
    \draw[->,thick] (x4.-112) .. controls (4.5,-0.7) and (5.5,-0.7) .. (x4.-68);
    \path[thick,->,draw,black]
    (y1) edge (x1)
    (y2) edge (x1)
    (y3) edge (x1)
    (y4) edge (x1)
    (y5) edge (x1)
    (y6) edge (x1)
    (y7) edge (x1)
    (y8) edge (x1)
    (y9) edge (x1)
    (y10) edge (x1)
    (y11) edge (x1)
    (y12) edge (x1)
    ;
  \end{tikzpicture}
  \caption{\label{fig:dense_f}
  The function of Lemma \ref{lem:k_non_sources} for $n = 4$ and $q =k= 2$.}
\end{figure}
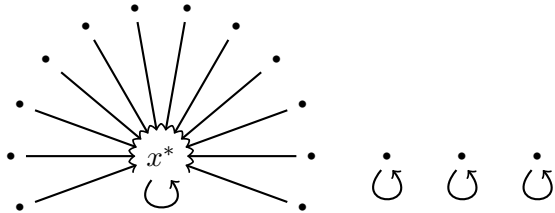

\begin{proof}
  Let  $X$ be a subset of $\Q{q}^n$ of size $q^k$ and let $x^*\in X$. Let $f
  \in F(n,q)$ such that for all $x \in \Q{q}^n$, $f(x) = x$ if $x \in X$ and
  $f(x)=x^*$ otherwise. See Figure \ref{fig:dense_f} for an illustration. The
  images of $f$ are thus exactly $X$. We will prove that $G(f)$ has at least
  $k$ non-sources, each of strict in-degree at least $n-k$. Since $f$ is
  described up to isomorphism, this proves the lemma. 

  \medskip
  We first prove that every non-source of $G(f)$ has strict in-degree at least
  $n-k$. Let $Y$ be any non-empty strict subset of $X$, and let
  $Y^-=f^{-1}(Y)$. We have $0<|Y|<|X|=q^k$ thus $|Y|$ is not a multiple of
  $q^k$. We will show that $|Y^-|$ is not a multiple of $q^k$. Note that
  $|f^{-1}(x^*)| = q^n-q^k+1$ and $|f^{-1}(x)| = 1$ for all $x\in
  X\setminus\{x^*\}$. Thus, if $x^*\notin Y$ then $|Y^-|=|Y|$ is not a
  multiple of $q^k$. If $x^*\in Y$ then
  $|Y^-|=|f^{-1}(x^*)|+|Y|-1=q^n-q^k+|Y|$, and since $|Y|$ is not a multiple of
  $q^k$ we deduce that $|Y^-|$ is not a multiple of $q^k$. Consequently, by
  Lemma~\ref{lem:sum_images}, every non-source of $G(f)$ has strict in-degree
  at least $n-k$. 

  \medskip
  Let $\ell$ be the number of non-sources of $G(f)$. It remains to prove that
  $\ell\geq k$. Let $I$ be the set of sources of $G(f)$. For each $i\in I$,
  $f_i$ is a constant function which always returns some member $a_i\in \Q{q}$.
  Let $Z$ be the set of $x\in \Q{q}^n$ with $x_i=a_i$ for all $i\in I$. We have
  $X\subseteq Z$, so $q^k=|X|\leq|Z|=q^{n-|I|}=q^\ell$. Thus, $\ell\geq k$ as
  desired. 
\end{proof}

The proof of Theorem~\ref{thm:dense} is now immediate: by Lemma
\ref{lem:k_non_sources} applied with $k=\lfloor n/2 \rfloor$, there is $f \in
F(n,q)$ such that every digraph in $\G(f)$ has at least $\lfloor n/2 \rfloor$
vertices of in-degree at least $n-\lfloor n/2 \rfloor$, and thus at least
$\lfloor n/2 \rfloor(n-\lfloor n/2 \rfloor)=\lfloor n^2/4 \rfloor$ arcs.

\section{Concluding remarks}\label{sec:concluding}

\begin{itemize}
  \item
    Let us say that a digraph $G$ is \EM{$(n,q)$-universal} if $G$ is
    isomorphic to some digraph in $\G(f)$ for every $f\in F(n,q)$ with
    $f\neq\id,\cst$. Theorem~\ref{thm:complete} says that $K_n$ is
    $(n,q)$-universal if $n\geq 5$ or $q\geq 3$. Actually, for $n,q\geq 2$, any
    $(n,q)$-universal digraph is isomorphic to $K_n$. Indeed, suppose $G$ is
    $(n,q)$-universal. First, if $f\in F(n,q)$ is a bijection, then every
    digraph in $\G(f)$ has no source. Thus, $G$ has no source. Then, by
    Lemma~\ref{lem:k_non_sources} applied with $k=1$, there is $f\in F(n,q)$
    such that any digraph in $\G(f)$ has at least one non-source, and each
    non-source has strict in-degree $n-1$. So $G$ has this property, and since
    it has no source, we deduce that each vertex has strict in-degree $n-1$.
    Finally, if $f\in F(n,q)$ has at least $q^{n-1}+1$ fixed points, then we
    easily check that any digraph in $\G(f)$ has $n$ loops. So $G$ has $n$
    loops and we deduce that it is isomorphic to $K_n$. 
  \item
    We proved that there is no $f\in F(n,q)$ such that $\G(f)$ only contains
    $K_n$ when $n$ is large compared to $q$, since we proved that
    $\delta^-(n,q)\leq c_q q$ for some constant $c_q$ that only depends on $q$,
    and thus $\delta^-(n,q)<n$ for $n>c_q q$. In particular, for $n\geq 5$,
    there is no $f\in F(n,2)$ such that $\G(f)=\{K_n\}$ since we proved that
    $\delta^-(5,2)\leq 4$ and $\delta^-(n,2)\leq 5$. However, for every $q\geq
    3$, there is $f\in F(2,q)$ such that $\G(f)=\{K_2\}$ since we proved that
    $\delta^-(2,q)= 2$. Does this phenomenon hold for every $n$, that is, for
    every fixed $n$ and $q$ large enough with respect to $n$, is there $f\in
    F(n,q)$ such that $\G(f)=\{K_n\}$? An equivalent formulation is: does
    $\delta^-(n,q)$ necessarily increase with $q$? 
    Also, we
    proved that $2\leq\delta^-(n,2)\leq 5$, for every $n\geq 3$, and we may ask
    if $\delta^-(n,2)$ is a constant, and if so which one of the four possible
    values it takes.
  \item
    Let $d(n,q)$ be the maximum real $d\in [0,1]$ such that there is $f\in
    F(n,q)$ such that every digraph in $\G(f)$ has $dn^2$ arcs, that is, has
    density $d$. We prove that $d(n,q)\geq 1/4-1/n^2$ where $1/n^2$ tends to $0$
    as $n$ tends to infinity, \emph{i.e.}~is $o(1)$. Hence, this lower bound is independent of the
    alphabet size $q$, and we may ask, as above, whether $d(n,q)$ is necessarily
    increasing with $q$.
\end{itemize}

\paragraph{Acknowledgments}

This work was supported by 
the \emph{Young Researcher} project ANR-18-CE40-0002-01 ``FANs''.

\bibliographystyle{plain}
\bibliography{BIB}

\end{document}